\documentclass[]{article}
\usepackage{amsmath,amsthm,amssymb,graphicx,enumerate,natbib}

\newtheorem{theorem}{Theorem}[section]
\newtheorem{lemma}[theorem]{Lemma}

\newtheorem{assumption}[theorem]{Assumption}
\newtheorem{remark}[theorem]{Remark}
\newcommand{\defn}{:=}
\newcommand{\one}{{\bf 1}}
\def\bbr{\mathbb{R}}
\numberwithin{equation}{section}

\author{Souvik Ghosh and Gennady Samorodnitsky}
\date{\today}
\title{The effect of memory on functional large deviations of infinite
  moving average processes}

\begin{document}
\maketitle

\begin{abstract}
The large deviations of an infinite moving average process with
exponentially light tails are very similar to those of an
i.i.d. sequence as long as the coefficients decay fast enough. If they
do not, the large deviations change dramatically. We study this
phenomenon in the context of functional large, moderate and huge
deviation principles. 
\end{abstract}

{\bf Key words:}  large deviations, long range dependence, long
memory, moving average, rate function, speed function\\

\par 

{\bf AMS (2000) Subject Classification :}
60F10  (primary), 60G10, 62M10  (secondary)  

\thanks{Research partially supported by  NSA grant MSPF-05G-049, ARO
grant W911NF-07-1-0078 and NSF training grant ``Graduate and Postdoctoral
Training in Probability and Its Applications''  at Cornell University}

\section{Introduction} \label{sec:intro}

We consider a (doubly) infinite moving average process $(X_n)$ defined 
by 
 \begin{equation} \label{e:process}
 X_n\defn \sum_{i=-\infty}^\infty \phi_i\, Z_{n-i},n\in    \mathbb{Z}.
\end{equation} 
The innovations $\{Z_i,i\in \mathbb{Z}\}$ are assumed to be 
i.i.d. $\mathbb{R}^d$-valued 
light-tailed random variables with 0 mean  and covariance matrix
$\Sigma$. In this setup square summability of the coefficients $(\phi_i)$ 
 \begin{equation} \label{e:sq.summ}
\sum_{i=-\infty}^\infty
   \phi_i^2<\infty
\end{equation} 
is well known to be necessary and sufficient for convergence of the
series in \eqref{e:process}. We assume \eqref{e:sq.summ} throughout
the paper. 
Under these assumption $(X_n)$ is a well defined stationary process,
also known as a linear process; see \cite{brockwell:davis:1991}. It is
common to think of a linear process  as a short memory process when it satisfies 
the stronger condition of absolute summability of coefficients,
\begin{equation}\label{absum}
\sum_{n\in\mathbb{Z}}|\phi_i|<\infty. 
\end{equation}
One can easily  check that absolute
summability of coefficients implies absolute
summability of the covariances: 
\[
	\sum_{i=-\infty}^\infty|Cov(X_0,X_i)|<\infty
\] 
It is also easy to exhibit a broad class of
examples where \eqref{absum} fails and the covariances are not
summable.

Instead of covariances, we are interested in understanding how the
large deviations of a moving average process change as the
coefficients decay slower and slower. Information obtained in this way is
arguably more substantial than that obtained via covariances alone. 

We assume that the moment generating function of a generic noise
variable $Z_0$, is finite in a neighborhood of the origin. We denote  its log-moment generating function by $\Lambda(\lambda)\defn \log
E\big(\exp(\lambda\cdot   Z_0)\big)$, where  $x\cdot
y$ is  the scalar product of two  vectors, $x$ and $y$. For 
 a function $f:\mathbb{R}^d\rightarrow (-\infty,\infty]$, 
define the Fenchel-Legendre transform of $f$ by $f^* =
\sup_{\lambda\in \mathbb{R}^d}\bigl\{ \lambda\cdot x - f(x)\bigr\}$, 
and the set $\mathcal{F}_f\defn \{ x\in \mathbb{R}^d:
f(x)<\infty \}\subset  
\mathbb{R}^d$. The imposed assumption  $0\in 
\mathcal{F}_\Lambda^\circ$, the interior of $\mathcal{F}_\Lambda$, is
then the formal statement of our comment that the innovations $(Z_i)$
are light-tailed. Section 2.2 in \cite{dembo:zeitouni:1998} summarizes
the  properties of $\Lambda$ and $\Lambda^*$. 

We are interested in the large deviations of probability measures 
based on partial sums of a moving average process. Recall that a
sequence of probability measures 
$\{\mu_n\}$ on the Borel subsets of a topological space 
is said to satisfy the \emph{large deviation principle}, or LDP,
with speed $b_n$, and upper 
and lower rate function $I_u(\cdot)$ and $I_l(\cdot)$, respectively,
if for any Borel set $A$, 
\begin{equation}\label{ldpdef}
-\inf_{x\in A^{\circ}}I_l(x)\le \liminf_{n\rightarrow \infty
 }\frac{1}{b_n}\log \mu_n(A)\le \limsup_{n\rightarrow
 \infty}\frac{1}{b_n}\log \mu_n(A)\le -\inf_{x\in \bar{A}}I_u(x), 
\end{equation}
where $ A^{\circ}$ and $\bar{A}$ are, respectively, the interior and
closure of $A$. 
A rate function is a non-negative lower semi-continuous function, and
a good rate function is a rate function with compact level sets.  We
refer the reader to \cite{varadhan:1984}, \cite{deuschel:stroock:1989} or 
\cite{dembo:zeitouni:1998}  for a detailed treatment of  large
deviations. 

In many cases, the sequence of
measures $\{\mu_n\}$ is the sequence of the laws of the normalized 
partial sums $a_n^{-1}(X_1+\ldots +X_n)$, for some appropriate
normalizing sequence $(a_n)$. Large deviations can also be formulated
in function spaces, or in measure spaces. The normalizing sequence has
to grow faster than the rate of growth required to obtain a
non-degenerate weak limit theorem for the normalized partial
sums. There is, usually, a boundary for the rate of growth of the
normalizing sequence, that separates the ``proper large deviations'' 
from the so-called ``moderate deviations''. In the moderate deviations
regime the normalizing sequence $(a_n)$ grows slowly enough so as to
make the underlying weak limit felt, and Gaussian-like rate functions
appear. This effect disappears at the boundary, which corresponds to
the proper large deviations. Normalizing sequences that grow even
faster lead to the so-called ``huge deviations''. For the
i.i.d. sequencies $X_1,X_2,\ldots$ the proper large deviations regime
corresponds to the linear growth of the normalizing sequence. The same
remains true for certain short memory processes. We will soon see  that for certain long memory processes the natural boundary is
not the linear normalizing sequence. 

There exists rich literature on large deviation for moving average
processes, going back to 
\cite{donsker:varadhan:1985}. They considered  Gaussian moving
averages and proved LDP for  the random  measures 
$n^{-1}\sum_{i\le n}\delta_{X_i}$, under the assumption that the
spectral density of the process is
continuous. \cite{burton:dehling:1990} considered a general
one-dimensional moving average 
process with $\mathcal{F}_\Lambda=\mathbb{R}$, assuming that 
\eqref{absum} holds. They also assumed that 
 \begin{equation}\label{e:nonzero}
 \sum_{n\in\mathbb{Z}}\phi_i=1; 
\end{equation}
the only substantial part of the assumption being that the sum of the
coefficients in non-zero. In that case $\{\mu_n\}$, the laws of  
$n^{-1}S_n = n^{-1}(X_1+\ldots +X_n)$, 
satisfy  LDP with a good rate function $\Lambda^*(\cdot)$. 
The work of \cite{jiang:rao:wang:1995} handled the case 
of  \mbox{$\{Z_i,i\in \mathbb{Z}\}$}, taking values in a separable Banach
space. Still assuming  (\ref{absum}) and \eqref{e:nonzero}, they
proved that the sequence $\{\mu_n\}$ satisfies a large deviation 
lower bound with the good rate function $\Lambda^*(\cdot)$, and, under
an integrability assumption, a large
deviation upper bound also holds with a certain good rate function
$\Lambda^\#(\cdot)$. In a finite dimensional Euclidian
space, the integrability assumption is equivalent to \mbox{$0\in 
\mathcal{F}^\circ_\Lambda$,} and the upper rate function is given  by
\begin{equation}
\Lambda^\#(x) \defn \sup_{\lambda\in \Pi} \{ \lambda\cdot
x-\Lambda(\lambda)\},
\end{equation} 
where  $\Pi=\{ \lambda\in \mathbb{R}^d$:  there exists $N_\lambda$
such that $\sup_{n\ge N_\lambda,i\in\mathbb{Z}}\Lambda(\lambda
\phi_{i,n})<\infty\}$ with $\phi_{i,n}\defn
\phi_{i+1}+\cdots+\phi_{i+n}$. Observe that, if 
$\mathcal{F}_\Lambda=\mathbb{R}^d$, then $\Lambda^\#\equiv
\Lambda^*$. 

In their paper, \cite{djellout:guillin:2001} went back to the
one-dimensional case. They worked under the assumption that the
spectral density is continuous and non-vanishing at the origin. 
Assuming also that the 
noise variables have a  bounded support, they showed that the LDP of 
\cite{burton:dehling:1990} still holds, and also established a
moderate deviation principle. 

\cite{wu:2004} extended the results of \cite{djellout:guillin:2001} and proved a large 
deviation principle for the occupation measures of the moving
average processes. He worked in an arbitrary dimension $d\geq 1$, with the
 same assumption on the spectral density  but  replaced the
 assumption of the boundedness of the  
support of the noise variables with the strong integrability condition,
$E[\exp(\delta|Z_0|^2)]<\infty$, for some $\delta>0$.  It is worth noting that  an explicit rate function could be obtained only under the absolute summability
assumption (\ref{absum}). 

Further, \cite{jiang:rao:wang:1992} considered moderate deviations 
in one dimension  under the absolute summability of the coefficients,
and assuming that $0\in  \mathcal{F}^\circ_\Lambda$. Finally,
\cite{dong:xi-li:yang:2005} showed that, under the same summability and
integrability assumptions,
the moving average ``inherits'' its moderate deviations from the noise
variables even if the latter are not necessarily i.i.d. 

Our main goal in this paper is to understand what happens when the
absolute summability of the coefficients  (or a variation, like
existence of a spectral density which is non-zero and continuous at
the origin) 
fails. Specifically, we will assume a certain regular variation
property of the coefficients; see Section \ref{sec:FLDP}.  For
comparison, we also present parallel results for the case 
where the coefficients are summable (most of the results are new even
in this case). We will see that there is a significant difference
between large deviations in the case of absolutely summable
coefficients (which are very similar to the large deviations of an
i.i.d. sequence) and the situation we consider, where absolute
summability fails.  In this sense, there is  a justification for
viewing (\ref{absum}), or ``its neighbourhood'', as the short memory
range of coefficients for a moving average process. Correspondingly,
the complementary situation may be viewed as describing the long 
memory range of coefficients for a moving average process. A similar
phenomenon occurs in important applications to 
{\it ruin probabilities} and {\it long strange segments}; a discussion
will appear in a companion paper.

The main part of the paper is Section \ref{sec:FLDP}, where we discuss
functional large deviation principles for a moving average process in
both short and long 
memory settings. Certain lemmas required for the proofs in that
section  are postponed until Section \ref{sec:lemmas}.

\section{Functional large deviation principle} \label{sec:FLDP}
This section discusses the large, moderate and huge
deviation principles for the sample paths of the moving average process. 
Specifically, we study the step process  $\{Y_n\}$ 
\begin{equation}\label{stepjump} 
Y_n(t)=\frac{1}{a_n}\sum_{i=1}^{[nt]}X_i,t\in [0,1],
\end{equation} 
and its polygonal path counterpart
\begin{equation}\label{e:polyg}
\tilde Y_n(t)=\frac{1}{a_n} \sum_{i=1}^{[nt]}X_i
+\frac{1}{a_n}(nt-[nt])X_{[nt]+1},t\in[0,1]. 
\end{equation}
Here $(a_n)$ is an appropriate normalizing sequence. 
We will use the notation $\mu_n$ and $\tilde \mu_n$ to denote the laws of
$Y_n$ and $\tilde Y_n$, respectively, in the function space appropriate to the situation at
hand, equipped with the cylindrical $\sigma$-field. 

Various parts of the theorems in this section will work with several
topologies on 
the space $\mathcal{BV}$ of all $\mathbb{R}^d$-valued functions  of
bounded variation defined on the unit interval $[0,1]$. To ensure that the
space $\mathcal{BV}$ is a measurable set in the cyindrical
$\sigma$-field of all $\mathbb{R}^d$-valued functions on $[0,1]$, we
use  only rational partitions of $[0,1]$ when defining variation.  We
will use subscripts to denote the topology on the 
space. Specifically, the subscripts $S$, $P$ and $L$ will denote the
sup-norm topology, the topology of
pointwise convergence and, finally, the topology in which $f_n$
converges to $f$ if and only if $f_n$ converges to $f$ both pointwise
and in $L_p$ for all $p\in [1,\infty)$.

We call a function
$f:\mathbb{R}^d\rightarrow \mathbb{R}$  \emph{balanced
  regular varying} with exponent $\beta>0$, if there exists a
non-negative bounded function $\zeta_f$
defined on the unit sphere on $\mathbb{R}^d$ and a function
$\tau_f:[0,\infty)\rightarrow   [0,\infty)$ satisfying  
\begin{equation}\label{bregvar1}
\lim_{t\rightarrow \infty}\frac{\tau_f(tx)}{\tau_f(t)}=x^\beta
\end{equation} 
for all $x>0$ (i.e. $\tau_f$ is regularly varying with exponent
$\beta$)  such that for any $(\lambda_t)\subset \mathbb{R}^d$
converging to $\lambda$, with $|\lambda_t|=1$ for all $t$, we have
\begin{equation}\label{bregvar2} 
\lim_{t\rightarrow \infty}\frac{
    f(t\lambda_t)}{\tau_f(t)}=\zeta_f(\lambda).
\end{equation} 
We will typically omit the subscript $f$ if doing so is not likely to
cause confusion. 

The following assumption describes the short memory scenarios we
consider. In addition to the summability of the coefficients, the
different cases arise from the ``size'' of the normalizing constants
$(a_n)$ in \eqref{stepjump}, the resulting speed sequence $(b_n)$ and
the integrability assumptions on the noise variables. 
\begin{assumption} \label{ass:short}
All the scenarios below assume that 
\begin{equation}\label{cabsum}
\sum_{i\in \mathbb{Z}}|\phi_i|<\infty
  \mbox{ and }\sum_{i\in \mathbb{Z}}\phi_i=1.
\end{equation} 
\begin{enumerate}[$S1.$]
\item $a_n=n,0\in \mathcal{F}_\Lambda^\circ$ and $b_n=n$.
\item $a_n=n,\mathcal{F}_\Lambda=\mathbb{R}^d$ and $b_n=n$.
\item $a_n/\sqrt{n}\rightarrow \infty, \ a_n/n\rightarrow 0, \ 0\in
  \mathcal{F}_\Lambda^\circ$ and $b_n=a_n^2/n$. 
\item $a_n/n\rightarrow \infty$, $\Lambda(\cdot)$ is balanced regular
  varying with exponent $\beta>1$ and $b_n=n\tau(\gamma_n)$, where  
\begin{equation}\label{cn} \gamma_n=\sup \{x:\tau(x)/x\le
  a_n/n\}.\end{equation} 
\end{enumerate}
\end{assumption}

Next, we introduce a new notation required to state our first result. 
For $i\in \mathbb{Z}$ and $n\ge 1$ we set 
$\phi_{i,n}\defn \phi_{i+1}+\cdots+\phi_{i+n}$. Also for $k\ge
1$ and $0<t_1<\cdots<t_k\le 1$, a subset $\Pi_{t_1,\ldots,t_k}\subset
(\mathbb{R}^d)^k$ is defined by 
$$
\Pi_{t_1,\ldots,t_k}\defn
\Big\{\underline{\lambda}=(\lambda_1,\ldots,\lambda_k) \in
(\mathcal{F}_\Lambda)^k: \ \Lambda \ \text{is continuous on
  $\mathcal{F}_\Lambda$ at each $\lambda_j$,} 
$$
\begin{equation}\label{pi}
\text{and for   some $N\ge 1$,} \ \sup_{n\ge
  N,\, j\in\mathbb{Z}}
\Lambda\Big(\sum_{i=1}^k\lambda_i\phi_{j+[nt_i], 
  [nt_i]-[nt_{i-1}]} \Big)<\infty \Big\}.  
\end{equation} 

We view the next theorem as describing the sample path large
deviations of (the partial sums of) a moving average process in the
short memory case. The long memory counterpart is theorem \ref{thmlm}
below. 

\begin{theorem}\label{thmsm}
\begin{enumerate}[(i)]
\item If $S1$ holds, then $\{\mu_n\}$ satisfy in $\mathcal{BV}_L$, LDP
  with speed $b_n\equiv n$, good upper rate function 
\begin{equation}\label{upperratefn}
G^{sl}(f)=\sup_{k\geq 1, \, t_1,\ldots,t_k}
\Big\{\sup_{\underline{\lambda}\in \Pi_{t_1,\ldots,t_k} }
\sum_{i=1}^{k}\Big\{\lambda_i\cdot\big(f(t_i)
-f(t_{i-1})\big)-(t_i-t_{i-1})\Lambda(\lambda_i) \Big\} \Big\} 
\end{equation} 
if $f(0)=0$ and $G^{sl}(f)=\infty$ otherwise, and with good lower
rate function  
\[
H^{sl}(f)=\left\{\begin{array}{ccl}  \int\limits_0^1
\Lambda^*(f^\prime(t))dt & if & f\in \mathcal{AC}, f(0)=0\\ \infty &
 & otherwise,\\ \end{array}\right.
\] 
where $\mathcal{AC}$ is the set of all absolutely continuous
functions, and $f^\prime$ is the coordinate-wise derivative of $f$.  
\item  If $S2$ holds, then $H^{sl}\equiv G^{sl}$ and $\{\mu_n\}$
  satisfy LDP in $\mathcal{BV}_S$, with speed $b_n\equiv n$ and good
  rate function   $H^{sl}(\cdot)$. 
\item Under assumption $S3$, $\{\mu_n\}$ satisfy in $\mathcal{BV}_S$,
  LDP with speed $b_n$ and good rate function   
\[
H^{sm}(f)=\left\{\begin{array}{ccl}  \int\limits_0^1
\frac{1}{2}f^\prime(t)\cdot \Sigma^{-1} f^\prime(t)dt & if & f\in
\mathcal{AC}, f(0)=0\\ \infty &  & otherwise.\\ \end{array}\right.
\]
Here   $\Sigma$ is the covariance matrix of $Z_0$, and we understand
$a\cdot \Sigma^{-1}a$ to mean $\infty$ if $a\in K_\Sigma\defn \{x\in \mathbb{R}^d-\{0\}:\Sigma x=0\}$. 
\item Under assumption $S4$, $\{\mu_n\}$ satisfy in $\mathcal{BV}_S$,
  LDP with   speed $b_n$ and good rate function
  \[
H^{sh}(f)=\left\{\begin{array}{ccl}  \int\limits_0^1
  (\Lambda^{h})^*(f^\prime(t))dt   & if & f\in \mathcal{AC}, f(0)=0\\
  \infty &  & otherwise,\\ \end{array}\right..
\] 
where $\Lambda^{h}(\lambda)=\zeta_\Lambda\Big(\frac{\lambda}{|\lambda|}\Big)
|\lambda|^\beta$ for $\lambda\in \mathbb{R}^d$ (defined as zero for
$\lambda=0$). 
\end{enumerate}
\end{theorem}

A comparison with the LDP for i.i.d. sequences (see
\cite{mogulski:1976} or theorem 5.1.2 in \cite{dembo:zeitouni:1998})
reveals that the rate function stays the same as long as the
coefficients in the moving average process stay summable. 

We also note that 
an application of the contraction principle gives, under scenario S1,
a marginal LDP  for the law of $n^{-1}S_n$ in   $\mathbb{R}^d$  with
speed $n$, upper rate function  $G^{sl}_1(x)=\sup_{\lambda\in
  \Pi_1}\Big\{\lambda\cdot x   -\Lambda(\lambda) \Big\}$, 
and lower rate function   $\Lambda^*(\cdot)$, recovering the statement
of theorem 1 in   \cite{jiang:rao:wang:1995} in the finite-dimensional
case.


Next, we consider what happens when
the absolute summability fails, in a ``major way''. We will assume
that the coefficients are balanced  regular varying with an
appropriate exponent. The following assumption is parallel to
assumption \ref{ass:short} in the present case, dealing, once again,
with the various cases that may arise. 

\begin{assumption}\label{regularvarying} \label{ass:long}
All the scenarios assume that the coefficients $\{\phi_i\}$ are
  balanced   regular varying with exponent $-\alpha,1/2<\alpha\le 1$
  and   $\sum\limits_{i=-\infty}^\infty|\phi_i|=\infty$. Specifically,
  there is   $\psi:[0,\infty)\rightarrow [0,\infty)$ and $0\le p\le
  1$, such that   
\begin{equation}\label{coeff}
  \left.\begin{array}{c}\lim\limits_{t\rightarrow
  \infty}\frac{\psi(tx)}{\psi(t)}= x^{-\alpha}, \ \text{for all $x>0$} \\   
 \lim\limits_{n\rightarrow\infty} \frac{\phi_n}{\psi(n) } =p \mbox{
  and } \lim\limits_{n\rightarrow \infty}  \frac{\phi_{-n}}{\psi(n)
  }=q\defn 1-p.\end{array} \right\} \end{equation} 
Let $\Psi_n\defn \sum_{1\le i\le n}\psi(i)$. 
\begin{enumerate}[$R1.$]
\item $a_n=n\Psi_n,0\in   \mathcal{F}_\Lambda^\circ$ and $b_n=n$.
\item $a_n=n\Psi_n,\mathcal{F}_\Lambda=\mathbb{R}^d$ and $b_n=n$.
\item $a_n/\sqrt{n}\Psi_n\rightarrow \infty, a_n/(n\Psi_n)\rightarrow
  0, 0\in \mathcal{F}_\Lambda^\circ$ and $b_n=a_n^2/(n\Psi_n^2)$. 
\item $a_n/(n\Psi_n)\rightarrow \infty,$ $\Lambda(\cdot)$ is balanced
  regular varying with exponent $\beta>1$ and
  $b_n=n\tau(\Psi_n\gamma_n)$, where  
\begin{equation}\label{cnp} \gamma_n=\sup \{x:\tau(\Psi_nx)/x\le
  a_n/n\}.\end{equation} 
\end{enumerate} 
\end{assumption}


Similar to \eqref{pi} we define 
$$
\Pi^{\alpha}_{t_1,\ldots,t_k} \defn
\Big\{\underline{\lambda}=(\lambda_1,\ldots,\lambda_k): (p\wedge
q)\lambda_i \in \mathcal{F}_\Lambda^\circ, \, i=1,\ldots, k,\ \text{and}
$$
\begin{equation}\label{pir}
\text{for   some $N=1,2,\ldots$} \ \sup_{n\ge
  N,\, j\in\mathbb{Z}}
\Lambda\Big(\frac{1}{\Psi_n}\sum_{i=1}^k\lambda_i\phi_{j+[nt_i], 
  [nt_i]-[nt_{i-1}]} \Big)<\infty \Big\}  
\end{equation} 
for $1/2<\alpha<1$, while for $\alpha=1$, we define 
$$
\Pi^{1}_{t_1,\ldots,t_k} \defn
\Big\{\underline{\lambda}= (\lambda_1,\ldots,\lambda_k) \in
(\mathcal{F}_\Lambda)^k: \ \Lambda \ \text{is continuous on
  $\mathcal{F}_\Lambda$ at each $\lambda_j$}
$$
\begin{equation}\label{pir.1}
\text{ and for   some $N=1,2,\ldots$} \ \sup_{n\ge
  N,\, j\in\mathbb{Z}}
\Lambda\Big(\frac{1}{\Psi_n}\sum_{i=1}^k\lambda_i\phi_{j+[nt_i], 
  [nt_i]-[nt_{i-1}]} \Big)<\infty \Big\}  
\end{equation} 

Also for $1/2<\alpha<1$, any $ k\ge 1$, $0< t_1\le \cdots\le
t_k\le 1$, and $\underline{\lambda}=(\lambda_1,\ldots,\lambda_k)\in
(\mathbb{R}^d)^k$ let 
\begin{equation} \label{e:h}
h_{t_1,\ldots,t_k}(x;\underline{\lambda}) \defn
(1-\alpha)\sum\limits_{i=1}^k\lambda_i
\int\limits_{x+  t_{i-1}}^{x+t_i} |y|^{-\alpha} 
( pI_{[y\ge 0]}+qI_{[y<0]})dy.
\end{equation}

For any $\mathbb{R}^d$-valued convex function $\Gamma$, any function
$\varphi\in L_1[0,1]$ and $1/2<\alpha<1$ we define , 
\begin{equation} \label{e:transf.alpha}
\Gamma_\alpha^\ast(\varphi) = \sup_{\psi\in L_\infty[0,1]} \bigg\{
\int_0^1 \psi(t)\cdot \varphi(t)\, dt 
\end{equation}
$$
- \int_{-\infty}^\infty 
\Gamma\left( \int_0^1 \psi(t)(1-\alpha)|x+t|^{-\alpha}\Bigl[ 
pI_{[x+t\ge 0]}+qI_{[x+t<0]}\Bigr]\, dt\right)dx\biggr\}\,,
$$
whereas for $\alpha=1$ we put
\begin{equation} \label{e:transf.1}
\Gamma_1^\ast(\varphi) = \int_0^1 \Gamma^\ast(\varphi(t))\, dt\,.
\end{equation}

We view the following result as describing the large deviations of
moving averages in the long memory case.

\begin{theorem}\label{thmlm}
\begin{enumerate}[(i)]
\item If $R1$ holds, then $\{\mu_n\}$ satisfy in $\mathcal{BV}_L$, LDP
  with speed $b_n=n$, good upper rate function 
\begin{equation} \label{e:upperratefn.1}
G^{rl}(f)=\sup_{k\ge 1, t_1,\ldots,t_k}
\Big\{\sup_{\underline{\lambda}\in \Pi^\alpha_{t_1,\ldots,t_k} }
\sum_{i=1}^{k}\lambda_i\cdot\big(f(t_i)
-f(t_{i-1})\big)-\Lambda^{rl}_{t_1,\ldots,t_k}(\lambda_1,\ldots,\lambda_k)
\Big\} 
\end{equation} 
if $f(0)=0$ and $G^{rl}(f)=\infty$ otherwise, where
\begin{eqnarray}\label{e:lambda.rl}
 \Lambda^{rl}_{t_1,\cdots,t_k}(\lambda_1,\cdots,\lambda_k)
& \defn & \left\{ \begin{array}{lcl}\int\limits_{-\infty}^\infty
   \Lambda\Big(h_{t_1,\ldots,t_k}(x;\underline{\lambda}) \Big) dx & if
   & \alpha<1\\ \sum\limits_{i=1}^{k}(t_i-t_{i-1})\Lambda(\lambda_i) &
   if & \alpha=1,\end{array}\right. 
\end{eqnarray}
and  good lower  rate function 
\[
H^{rl}(f)= \left\{\begin{array}{ccl}
\Lambda_\alpha^\ast(f^\prime) & if & f\in \mathcal{AC}, f(0)=0\\ \infty &
 & otherwise.\\ \end{array}\right.
\] 


\item  If $R2$ holds, then $H^{rl}\equiv G^{rl}$ and $\{\mu_n\}$
  satisfy LDP in $\mathcal{BV}_S$, with speed $b_n=n$ and good rate
  function   $H^{rl}(\cdot)$. 
\item Under assumption $R3$, $\{\mu_n\}$ satisfy in $\mathcal{BV}_S$,
  LDP with speed $b_n$ and good rate function   
\[H^{rm}(f)=\left\{\begin{array}{ccl}
(G_\Sigma)_\alpha^\ast(f^\prime) & if & f\in \mathcal{AC}, f(0)=0\\ \infty &
 & otherwise,\\ \end{array}\right.
\] 
where $G_\Sigma(\lambda) = \frac12 \lambda\cdot\Sigma\lambda$,
$\lambda\in \mathbb{R}^d$. 

equation

\item Under assumption $R4$, $\{\mu_n\}$ satisfy in $\mathcal{BV}_S$,
  LDP with speed $b_n$ and good rate function
\[
H^{rh}(f)=\left\{\begin{array}{ccl}
(\Lambda^{h})_\alpha^\ast(f^\prime) & if & f\in \mathcal{AC},
f(0)=0\\ \infty &  & otherwise,\\ \end{array}\right.
\] 
with $\Lambda^{h}$ as in theorem \ref{thmsm}. 


\end{enumerate}
\end{theorem}

We note that a functional LDP under the assumption $R2$, but 
for a non-stationary fractional ARIMA model was obtained by
\cite{barbe:broniatowski:1998}.

\begin{remark}\label{1dimldpr}
\emph{ The proof of theorem \ref{thmlm} below shows that, under the
  assumption $R1$, the laws of $(n\Psi_n)^{-1}S_n$ satisfy LDP
  with speed $n$, good lower rate function $\Lambda_1^{rl*}(\cdot)$
  and good upper rate function  $G^{rl}_1(x)\defn\sup_{\lambda\in
  \Pi^{\alpha}_1}\big\{\lambda \cdot x-\Lambda^{rl}_1(\lambda) \big\}$.
If  $R2$ holds, then
  $\Pi^{\alpha}_1=\mathbb{R}^d$ and $G^{rl}_1\equiv (\Lambda^{rl}_1)^*$.} 
\end{remark}

\begin{remark}\label{1dimmdpr}
\emph{It is interesting to note that under the assumption  $R3$ it is
  possible to choose $a_n=n$, and, hence, compare the large deviations
  of the sample means of moving average processes with summable and
  non-summable coefficients. We see that the sample means of moving
  average processes with summable coefficients satisfy LDP with speed
  $b_n=n$, while the sample means of moving
  average processes with non-suumable coefficients (under assumption
  $R3$) satisfy LDP 
  with speed $b_n=n/\Psi_n^2$, which is regular varying with exponent
  $2\alpha-1$. The markedly slower speed function in the latter case
  (even for $\alpha=1$ one has $b_n=nL(n)$, with a slowly varying
  function $L(\cdot)$ converging to zero) demonstrates a phase transition
  occurring here.} 
\end{remark}

\begin{remark}\label{rk:gaussian.rate.long}
{\rm Lemma \ref{l:gaussian.rate.long} at the end of this section
  describes certain properties of the rate function
  $(G_\Sigma)_\alpha^\ast$, which is, clearly, also the rate function
  in {\it all} scenarios in the Gaussian case.} 
\end{remark}

The proofs of theorems \ref{thmsm} and \ref{thmlm} rely on lemmas
 appearing in section 3.

\begin{proof}[Proof of theorem \ref{thmsm}]
(ii),\, (iii) and (iv):\  Let $\mathcal{X}$ be the set of all
  $\mathbb{R}^d$-valued functions defined on the unit interval
  $[0,1]$ and let $\mathcal{X}^o$ be the subset of $\mathcal{X}$, of
  functions which start at the origin. Define $J$ as the 
  collection  of  all ordered finite subsets of $(0,1]$ with a partial
    order defined by inclusion. For any $j=\{0<t_1<\ldots<t_{|j|}\le
    1\}$ define the projection $p_j:\mathcal{X}^o\rightarrow
    \mathcal{Y}_j$ as  $p_j(f)=(f(t_1),\ldots,f(t_{|j|}))$,  
    $f\in \mathcal{X}^o$. So $\mathcal{Y}_j$ can be identified with the
    space $(\mathbb{R}^d)^{|j|}$ and the projective limit of
    $\mathcal{Y}_j$ over $j\in J$  can be identified with
    $\mathcal{X}^o$ equipped with the topology of pointwise
  convergence. Note that $\mu_n\circ p_j^{-1}$ is the law of 
    \[Y_n^j=(Y_n(t_1),\ldots,Y_n(t_{|j|}))\] and let
    \begin{equation}\label{vn}V_n=\big(Y_n(t_1),Y_n(t_2)-Y_n(t_1),
\cdots,Y_n(t_{|j|})- Y_n(t_{|j|-1})\big).\end{equation} 
 By lemma  \ref{limits} we see  that for
    any $\underline{\lambda}=(\lambda_1,\ldots,\lambda_{|j|})\in
    (\mathbb{R}^d)^{|j|}$ 
$$
\lim_{n\rightarrow\infty}\frac{1}{b_n}\log
  E\big(\exp\big[b_n\underline{\lambda}\cdot V_n \big]\big)
=\lim_{n\rightarrow\infty}\frac{1}{b_n}\log
  E\exp\Big[\frac{b_n}{a_n}\sum_{i=1}^{|j|}\lambda_i\cdot
  \Big(\sum_{k=[nt_{i-1}]+1}^{[nt_{i}]} X_k \Big)\Big]
$$
$$
 = 
  \lim_{n\rightarrow\infty}\frac{1}{b_n}\sum_{l=-\infty}^{\infty}\Lambda
  \Big(\frac{b_n}{a_n}\sum_{i=1}^{|j|}\lambda_i\phi_{l+[nt_{i-1}],[nt_{i}]-
  [nt_{i-1}]}\Big)
$$
$$
 = 
  \sum_{i=1}^{|j|}(t_i-t_{i-1})\Lambda^v(\lambda_i)
 \defn  \Lambda^v_{t_1,\ldots,t_{|j|}}(\underline{\lambda}),
$$
where $t_0=0$ and for any $\lambda\in \mathbb{R}^d$,
\[ \Lambda^v(\lambda)=\left\{\begin{array}{lll} \Lambda(\lambda) &
\text{in part (ii),}\\ 
\frac{1}{2}\lambda\cdot \Sigma \lambda & \text{in part (iii),}\\
\zeta\Big(\frac{\lambda}{|\lambda|}\Big)|\lambda|^\beta & \text{in
  part (iv).}\\ 
\end{array} \right.
\] 
By the Gartner-Ellis theorem, the laws of $(V_n)$ satisfy LDP with speed
$b_n$ and good  rate function
\[
\Lambda^{v*}_{t_1,\cdots,t_{|j|}}(w_1,\ldots,w_{|j|})=\sum_{i=1}^{|j|}
(t_i-t_{i-1})\Lambda^{v*}\Big(\frac{w_i}{t_i-t_{i-1}}\Big), 
\]
where $(w_1,\ldots,w_{|j|})\in (\mathbb{R}^d)^{|j|}$. The map
$V_n\mapsto Y_n^j$ from $(\mathbb{R}^d)^{|j|}$ onto itself is one to
one and continuous. Hence the contraction principle tells us that 
$\{\mu_n\circ p_j^{-1}\}$
satisfy LDP in $(\mathbb{R}^d)^{|j|}$ with   good rate function  
\begin{equation}\label{Hsfindim}
H^{v}_{t_1,\ldots,t_{|j|}}(y_1,\ldots,y_{|j|})\defn
\sum_{i=1}^{|j|}(t_i-t_{i-1})\Lambda^{v*}\Big(\frac{y_i-y_{i-1}}
    {t_i-t_{i-1}}\Big),\end{equation}  
where we take $y_0=0$. By lemma \ref{l:expon.equiv}, the same holds for the measures
\mbox{$\{\tilde \mu_n\circ p_j^{-1}\}$}. Proceeding as in lemma 5.1.6 in
\cite{dembo:zeitouni:1998} this implies  
that the measures $\{\tilde \mu_n\}$ satisfy
LDP in the space $\mathcal{X}^o$ equipped with the topology of pointwise
convergence, with speed $b_n$ and the rate function
described in the appropriate part of the theorem. As $\mathcal{X}^o$ is a closed subset of $\mathcal{X}$, the same holds for $\{\tilde\mu_n\}$ in $\mathcal{X}$ and the rate function is infinite outside $\mathcal{X}^o$.
 Since $\tilde \mu_n(\mathcal{BV})=1$ for all $n\ge 1$ and
the 3 rate functions in parts (ii), (iii) and (iv) of the theorem are
infinite outside of $\mathcal{BV}$, we conclude that $\{\tilde \mu_n\}$
satisfy 
LDP in $\mathcal{BV}_P$ with the same rate function. The sup-norm
topology on $\mathcal{BV}$ is stronger than that of pointwise
convergence  and  by lemma \ref{expttns}, $\{\tilde \mu_n\}$ is
exponentially tight in $\mathcal{BV}_S$. So by corollary 4.2.6 in 
\cite{dembo:zeitouni:1998}, $\{\tilde \mu_n\}$ satisfy LDP in
$\mathcal{BV}_S$ with speed $b_n$ and good rate function
$H^v(\cdot)$. Finally, applying lemma \ref{l:expon.equiv} once again,
we conclude that the same is true for the sequence $\{ \mu_n\}$.

(i): We use the above notation. It follows from lemma  \ref{limits}
that for any partition $j$ of $(0,1]$ and
  $\underline{\lambda}=(\lambda_1,\ldots,\lambda_{|j|})\in
  (\mathbb{R}^d)^{|j|}$, 
\[
\limsup_{n\rightarrow \infty}\frac{1}{n} \log E\big[
  \exp\big(n\underline{\lambda}\cdot V_n\big) \big] \le
\chi(\underline{\lambda}),
\] 
where 
\[
\chi(\underline{\lambda})= \left\{ \begin{array}{ccc}
  \sum\limits_{i=1}^{|j|}(t_i-t_{i-1})\Lambda(\lambda_i)  & \text{if} \ \,
  \underline{\lambda}\in \Pi_{t_1,\ldots,t_{|j|}}\\ 
  \infty &    \text{otherwise}.\\ 
 \end{array} \right.
 \] 
The law of $V_n$ is exponentially tight since by
\cite{jiang:rao:wang:1995} the law of $Y_n(t_i)-Y_n(t_{i-1})$ is
exponentially tight in $\mathbb{R}^d$ for every $1\le i\le |j|$. Thus
by theorem 2.1 of 
\cite{deacosta:1985} the laws of $(V_n)$ satisfy a LD upper
bound with speed $n$ and rate function
\[
\sup_{\underline{\lambda}\in\Pi_{t_1,\ldots,t_{|j|}}}\left\{
\underline{\lambda}\cdot
\underline{w}-\sum_{i=1}^{|j|}(t_i-t_{i-1})\Lambda(\lambda_i)
\right\},
\]
which is, clearly, good. 
Therefore, the laws of $(Y_n(t_1),\ldots,Y_n(t_{|j|}))$ satisfy a LD
upper bound with speed $n$ and good rate function 
\begin{equation} 
G^{sl}_{t_1,\ldots,t_{|j|}}(\underline{y})\defn
\sup_{\underline{\lambda}\in\Pi_{t_1,\ldots,t_{|j|}}}\Big\{
\sum_{i=1}^{|j|}\lambda_i\cdot
(y_i-y_{i-1})-\sum_{i=1}^{|j|}(t_i-t_{i-1})\Lambda(\lambda_i)\Big\}.
\end{equation} 
Using the upper bound part of the Dawson-Gartner theorem, we see that
$\{\mu_n\}$ satisfy LD upper bound in $\mathcal{X}_P^o$ with speed $n$
and good rate rate function 
\[ 
G^{sl}(f)=\sup_{j\in J}
G^{sl}_{t_1,\ldots,t_{|j|}}\big(f(t_1),\cdots,f(t_{|j|})\big)
 \]  
and, as before, the same holds in $\mathcal{X}_P$ as well.

Next we prove that $(Y_n(t_1),\ldots,Y_n(t_{|j|}))$ satisfy a LD lower
bound with speed $n$ and rate function
$H^{v}_{t_1,\ldots,t_{|j|}}(\cdot)$ defined in (\ref{Hsfindim}) for
part (ii).  Let  
$$
V_n^\prime= \frac{1}{n}\Big(\sum_{|i|\le 2n}\phi_{i,[nt_1]}Z_{-i},
\sum_{|i|\le 2n}\phi_{i+[nt_1],[nt_2]-[nt_1]}Z_{-i},\cdots,
$$
$$
\sum_{|i|\le
  2n}\phi_{i+[nt_{|j|-1}],[nt_{|j|}]-[nt_{|j|-1}]}Z_{-i}\Big)
$$
and observe that the laws of $(V_n)$ and of $(V_n^\prime)$ are
exponentially equivalent. 


For $k>0$ large enough so that $p_k \defn P(|Z_0|\le k)>0$ we 
let $\mu_k = E\bigl(Z_0\bigl|\, |Z_0|\leq k\bigr)$, and note
that $|\mu_k|\to 0$ as $k\to\infty$.

\pagebreak
\noindent Let
$$
V_n^{\prime,k}= \frac{1}{n}\Big(\sum_{|i|\le
  2n}\phi_{i,[nt_1]}(Z_{-i}-\mu_k) ,
\sum_{|i|\le 2n}\phi_{i+[nt_1],[nt_2]-[nt_1]}(Z_{-i}-\mu_k),\cdots,
$$
$$
\sum_{|i|\le
  2n}\phi_{i+[nt_{|j|-1}],[nt_{|j|}]-[nt_{|j|-1}]}(Z_{-i}-\mu_k))\Big) 
:= V_n^{\prime} -a_{n,k} \,,
$$
where $a_{n,k} = (b_1^{(n)}\mu_k,b_2^{(n)}\mu_k,\ldots,
b_{|j|}^{(n)}\mu_k)\in (R^d)^{|j|}$ with some $|b_i^{(n)}|\leq c$, a
constant independent of 
$i$ and $n$. We define a new probability measure  
\[
\nu_{n}^k(\cdot)=P\Big(V_n^{\prime,k} \in \cdot, |Z_i|\le k,\
\text{for all} \ |i|\le 2n \Big)p_k^{-(4n+1)}\,.
\] 
Note that for all $\underline{\lambda}\in (\mathbb{R}^d)^{|j|}$ by
(the proof of part (i) of) lemma \ref{limits}, 
\begin{eqnarray*}
& \lim\limits_{n\rightarrow\infty} \frac{1}{n}\log\Big\{ p_k^{-(4n+1)}E\Big[
    \exp  \big(n\underline{\lambda}\cdot V_n^\prime \big) I_{[|Z_i|\le
	k,\ |i|\le 2n]} \Big]\Big\}&\\ 
&  =
    \sum\limits_{l=1}^{|j|}(t_l-t_{l-1})\Bigl( L^k(\lambda_l)
    -\lambda_l\mu_k\Bigr) -t_{|j|}\log p_k,& 
\end{eqnarray*}
where $L^k(\lambda)\defn \log E\big[\exp(\lambda\cdot
      Z_0)I_{[|Z_0|\le k]} \big]$, and so for every $k\ge 1$,
      $\{\nu_n^k,n\ge 1\}$ satisfy LDP  with speed $n$ and good rate
      function   
\begin{eqnarray}\label{approxrtfn}
&\sup\limits_{\underline{\lambda}}\Big\{\underline{\lambda}\cdot
  \underline{x}-  \sum\limits_{l=1}^{|j|}(t_l-t_{l-1})\Bigl(
  L^k(\lambda_l) -\lambda_l\mu_k\Bigr)
  \Big\} +t_{|j|}\log p_k&\nonumber \\  
&=\sum\limits_{l=1}^{|j|}(t_l-t_{l-1})L^{k*}
\Big(\frac{x_l+t_{|j|}\mu_k}{t_l-t_{l-1}}\Big)   +t_{|j|}\log p_k\,. &   
\end{eqnarray}
Since for any open set $G$ 
$$
\liminf_{n\rightarrow \infty } \frac{1}{n} \log P(V_n^{\prime,k}\in G) 
\geq \liminf_{n\rightarrow \infty } \frac{1}{n} \log \nu_n^k(G) +
4\log p_k\,,
$$ 
we conclude that for any $x$ and $\epsilon>0$, for all $k$ large
enough,
$$
\liminf_{n\rightarrow \infty } \frac{1}{n}\log P(V_n^\prime\in
B(\underline{x},2\epsilon))\ge  \liminf_{n\rightarrow \infty }
\frac{1}{n} \log \nu_n^k(B(\underline{x},\epsilon))+4\log p_k\,,
$$ 
where $B(\underline{x},\epsilon)$ is an open ball centered at $x$ with
radius $\epsilon$. 

Now note that for every $\lambda\in \mathbb{R}^d$, $L^k(\lambda)$ is
increasing to $\Lambda(\lambda)$ with $k$. So by theorem B3 in
\cite{deacosta:1988},  there exists
$\{\underline{x}^k\}\subset (\mathbb{R}^d)^{|j|}$, such that
$\underline{x}^k\rightarrow \underline{x}$, and   
\[
\limsup_{k\rightarrow \infty } \sum\limits_{l=1}^{|j|}(t_l-t_{l-1})L^{k*}
\Big(\frac{x_l^k}{t_l-t_{l-1}}\Big) \leq
\sum\limits_{l=1}^{|j|}(t_l-t_{l-1})L^{*}
\Big(\frac{x_l}{t_l-t_{l-1}}\Big)\,.
\] 
Since $\underline{x}^k-t_{|j|}\underline{\mu}_k\in
  B(\underline{x},2\epsilon)$ for $k$ large, where $\underline{\mu}_k =
  (\mu_k, \ldots, \mu_k)\in (R^d)^{|j|}$, we conclude that 
\[ \liminf_{n\rightarrow \infty } \frac{1}{n}\log P(V_n^\prime\in
  B(\underline{x},\epsilon))\ge
  -\sum_{l=1}^{|j|}(t_l-t_{l-1})\Lambda^*\Big(
  \frac{x_l}{t_l-t_{l-1}}\Big) \,.
\] 
Furthermore, because the laws of $(V_n)$ and of $(V_n^\prime)$ are
exponentially equivalent, the same statement holds with $V_n$
replacing $V_n^\prime$. We have, therefore, established that the laws
of $(Y_n(t_1),\ldots,Y_n(t_{|j|}))$ 
  satisfy a LD lower bound with speed $n$ and good rate function
  $H^{v}_{t_1,\ldots,t_{|j|}}(\cdot)$ defined in (\ref{Hsfindim}) for
part (ii). By the lower bound part of the Dawson-G\"artner theorem, 
  $\{\mu_n\}$ satisfy a LD lower bound in $\mathcal{X}_P$ with speed $n$
  and rate function  $\sup_{j\in J}
  H^{v}_{t_1,\ldots,t_{|j|}}(f(t_1),$ $\ldots,f(t_{|j|}))$. This rate
  function is identical to $H^{sl}$. 

 Notice that the lower rate function  $H^{sl}$ is infinite outside of
  the space  $\cap_{p\in [1,\infty)} L_p[0,1]$, and by lemma
  \ref{l:spprt.BV}, the same is 
  true for the upper rate function $G^{sl}$ (we view $\cap_{p\in
  [1,\infty)} L_p[0,1]$ as a 
  measurable subset of $\mathcal{X}$ with respect to the universal
  completion of the cylindrical $\sigma$-field). We conclude that the
  measures $\{\mu_n\}$ satisfy a LD   lower bound in $\cap_{p\in
  [1,\infty)} L_p[0,1]$ with the topology of pointwise
  convergence. Since this topology is coarser than the $L$ topology,
  we can use lemma \ref{expttnsl1} to conclude that the LD
  upper bound  and the LD lower bound also hold in $\cap_{p\in
  [1,\infty)} L_p[0,1]$
  equipped with   $L$ topology. Finally, 
  the rate functions are also infnite outside of the space $\mathcal
  {BV}$,   and so the measures $\{\mu_n\}$ satisfy the LD  bounds in 
  $\mathcal {BV}$ equipped with $L$ topology. 
\end{proof}

\begin{proof}[Proof of theorem \ref{thmlm}]
The proof of parts (ii), (iii) and (iv) is identical to the proof of
the corresponding parts in theorem \ref{thmsm}, except that now lemma 
\ref{limitr} is used instead of lemma  \ref{limits}, and we
use lemma \ref{l:funct.rate.long} to identify the rate function. 

We now prove  part (i) of the theorem. We start by proving the finite
dimensional LDP for the laws of $V_n$ in (\ref{vn}). An inspection
of the proof of the corresponding statement on theorem \ref{thmsm}
shows that the only 
missing ingredient needed to obtain the upper bound part of this LDP 
is the exponential tightness of $Y_n(1)$ in $\mathbb{R}^d$. 
Notice that for $s>0$ and small $\lambda>0$ 
\[ 
P\Big( Y_n(1) \notin [-s,s]^d \Big)\le e^{-\lambda n s}
\sum_{l=1}^dE\Big( e^{\lambda  Y_n^{(l)}(1)}+e^{-\lambda  Y_n^{(l)}(1)}
\Big), 
\]
where $Y_n^{(l)}(1)$ is the $l$th coordinate of $Y_n(1)$. Since $0\in \mathcal{F}_\Lambda^o$, by part (i) of lemma
\ref{limitr} we see that  
\[ 
\lim_{s\rightarrow \infty}\limsup_{n\rightarrow \infty}
\frac{1}{n}\log P\Big( Y_n(1)\notin [-s,s]^d \Big)
=-\infty\,,
\]
which is the required exponential tightness. It follows  that the laws
of $(V_n)$ satisfy a LD upper bound with speed $n$ and rate function  
\[
\sup_{\underline{\lambda}\in \Pi^{rl}_{t_1,\ldots,t_{|j|}}}\left\{
\underline{\lambda}\cdot
\underline{w}-\Lambda^{rl}_{t_1,\ldots,t_{|j|}}
(\lambda_1,\ldots,\lambda_{|j|}) \right\}.
\]
Next we prove a LD lower bound for the laws of $(V_n)$.  The proof
in the case $\alpha=1$ follows the same steps as the corresponding
argument in theorem \ref{thmsm}, so we will concentrate on the case
$1/2<\alpha<1$. For $m\geq 1$ let 
\[
V_{n,m}^\prime= \frac{1}{n\Psi_n}\Big(\sum_{|i|\le
  mn}\phi_{i,[nt_1]}Z_{-i}, \sum_{|i|\le
  mn}\phi_{i+[nt_1],[nt_2]-[nt_1]}Z_{-i},\cdots,
\]
\[
\sum_{|i|\le
  mn}\phi_{i+[nt_{|j|-1}],[nt_{|j|}]-[nt_{|j|-1}]}Z_{-i}\Big)\,.
\] 
Observe that $V_n = V_{n,m}^\prime+R_{n,m}^\prime$ for some
$R_{n,m}^\prime$ independent of $V_{n,m}^\prime$ and such that
for every $m$, $R_{n,m}^\prime\to 0$ in probability as
$n\to\infty$. We conclude that for any
$\underline{x}=(x_1,\cdots,x_{|j|})\in (\mathbb{R}^d)^{|j|}$,
$\epsilon>0$,  and $n$ sufficiently large, one  has
\begin{equation}\label{approxr} 
P(V_n\in B(\underline{x},2\epsilon) )\ge \frac{1}{2}
P(V_{n,m}^\prime\in B(\underline{x},\epsilon) ) \,.
\end{equation}

For $k\geq 1$ we define $p_k$ and $\mu_k$ as in the proof of theorem
\ref{thmsm}, and once again we choose $k$ large enough so that
$p_k>0$. We also define
$$
V_{n,m}^{\prime,k}=\frac{1}{n\Psi_n}\Big(\sum_{|i|\le
  mn}\phi_{i,[nt_1]}\bigl( Z_{-i}-\mu_k\bigr), \sum_{|i|\le
  mn}\phi_{i+[nt_1],[nt_2]-[nt_1]}\bigl( Z_{-i}-\mu_k\bigr),\cdots,
$$
$$
\sum_{|i|\le
  mn}\phi_{i+[nt_{|j|-1}],[nt_{|j|}]-[nt_{|j|-1}]}\bigl(
Z_{-i}-\mu_k\bigr)\Big) := V_{n,m}^{\prime} -a_{n,k}^{(m)} \,,
$$
where $a_{n,k}^{(m)} = (b_1^{(n,m)}\mu_k,b_2^{(n,m)}\mu_k,\ldots,
b_{|j|}^{(n,m)}\mu_k)\in (R^d)^{|j|}$ with some $|b_i^{(n,m)}|\leq
c_m$, a constant independent of $i$ and $n$. 

Once again we define  a new probability measure by 
\[
\nu_{n}^{k,m}(\cdot)=P\Big(V_{n,m}^{\prime,k} \in \cdot, | Z_i|\le
k,\  \text{for all} \   |i|\le mn \Big)p_k^{-(2mn+1)}\,.
\] 
Note that for all $\underline{\lambda}\in (\mathbb{R}^d)^{|j|}$, by
(the proof of) lemma \ref{limitr}, 
\begin{eqnarray*}
&& \lim_{n\rightarrow\infty} \frac{1}{n}\log\Big\{
  p_k^{-(2mn+1)}E\Big[ \exp  \big(n\underline{\lambda}\cdot
    V_{n,m}^{\prime,k} \big) I_{[| Z_i|\le _k,\, |i|\le mn]}
    \Big]\Big\}\\ 
& = & \int\limits_{-m}^mL^k\left(
  (1-\alpha)\sum\limits_{i=1}^{|j|}\lambda_i\int\limits_{x+
  t_{i-1}}^{x+t_i} |y|^{-\alpha} ( pI_{[y\ge 0]}+qI_{[y<0]})dy\right)dx \\
& - & (1-\alpha)\sum_{l=1}^{|j|} \lambda_i\cdot\mu_k \int\limits_{-m}^m
\left( \int\limits_{x+
  t_{i-1}}^{x+t_i} |y|^{-\alpha} ( pI_{[y\ge
  0]}+qI_{[y<0]})dy\right)dx -2m\log p_k\\ 
& = &Q^{k,m}(\underline{\lambda}) -\mu_k\cdot R^m(\underline{\lambda}) 
-2m\log p_k \ \text{(say)}
\end{eqnarray*}
where $L^k(\lambda)= \log E\big[\exp(\lambda\cdot Z_0)I_{[|Z_0|\le
      k]} \big],$ as defined before. Therefore, for every $k\ge 1$, $\{\nu_n^{k,m},n\ge 1\}$
satisfy LDP  with speed $n$ and good rate function
$(Q^{k,m})^*(\underline{x}-\underline{c}_{k,m})+2m\log p_k$, where 
$\underline{c}_{k,m}=(c_1^{m}\mu_k,c_2^{m}\mu_k,\ldots,
      c_{|j|}^{m}\mu_k) \in (\mathbb{R}^d)^{|j|}$ with
$$
c_i^{m} = (1-\alpha) \int_{-m}^m \left( \int\limits_{x+
  t_{i-1}}^{x+t_i} |y|^{-\alpha} ( pI_{[y\ge
  0]}+qI_{[y<0]})dy\right)dx\,.
$$
Note that for every $\lambda\in \mathbb{R}^d$,
$L^k(\lambda)$ is 
increasing to $\Lambda(\lambda)$  and $Q^{k,m}(\underline{\lambda})$
is increasing to 
$$
\Lambda^{rl, m}_{t_1,\ldots,t_{|j|}}(\underline{\lambda}) = 
\int\limits_{-m}^m 
   \Lambda\Big(h_{t_1,\ldots,t_k}(x;\underline{\lambda}) \Big) dx
$$
with $k$. 

An application of  theorem B3 in \cite{deacosta:1988}
shows, as in the proof of theorem \ref{thmsm}, that for any ball
centered at $x$ with radius $\epsilon$ 
\[ 
\liminf_{n\rightarrow \infty } \frac{1}{n}\log P(V_{n,m}^\prime\in
B(\underline{x},\epsilon))\ge
-(\Lambda^{rl,m}_{t_1,\ldots,t_{|j|}})^*(\underline{x}).
\] 
Appealing to \eqref{approxr} gives us 
$$
\liminf_{n\rightarrow \infty } \frac{1}{n}\log P(V_{n}\in
B(\underline{x},2\epsilon))\ge
-(\Lambda^{rl,m}_{t_1,\ldots,t_{|j|}})^*(\underline{x})
$$
for all $m\geq 1$. We now apply the above argument once again: for
every $\lambda\in \mathbb{R}^d$, $\Lambda^{rl,
  m}_{t_1,\ldots,t_{|j|}}(\underline{\lambda})$ increases to  $\Lambda^{rl
  }_{t_1,\ldots,t_{|j|}}(\underline{\lambda})$, and yet another
appeal to theorem B3 in \cite{deacosta:1988} gives us the desired 
LD lower bound for the laws of $(V_n)$ in the case $1/2<\alpha<1$. 

Continuing as in the proof of theorem \ref{thmsm} we conclude that 
 $\{\mu_n\}$ satisfy a LD lower bound in $\mathcal{X}_P$ with speed $n$
  and rate function  $\sup_{j\in J}
 (\Lambda^{rl}_{t_1,\ldots,t_{|j|}})^*(f(t_1),$ $f(t_2)-f(t_1),\ldots,f(t_{|j|})-f(t_{|j|-1}))$.  
By lemma \ref{l:funct.rate.long} this is equal to $H^{rl}(f)$ in the
 case $1/2<\alpha<1$, and in the case $\alpha=1$ the corresponding
 statement is the same as in theorem \ref{thmsm}. The fact that the LD
 lower bound  holds also in $\mathcal{BV}_L$ follows in the same way
 as in theorem \ref{thmsm}. This completes the proof. 
\end{proof}

The next lemma discusses some properties of the rate function
$(G_\Sigma)_\alpha^\ast$ in theorem \ref{thmlm}. For $0<\theta<1$, let 
$$
H_\theta = \left\{ \psi:\, [0,1]\to\bbr^d, \ \text{measurable,
  and} \ \  \int_0^1\int_0^1 \frac{|\psi(t)|
|\psi(s)|}{|t-s|^\theta} dt\, ds<\infty\right\}\,.
$$
If $\Sigma$ is a nonnegative definite matrix, we define an inner
product on $H_\theta$ by 
$$
(\psi_1,\psi_2)_\Sigma = \int_0^1\int_0^1 \frac{\psi_1(t)\cdot\Sigma
 \psi_2(s)}{|t-s|^\theta} dt\, ds\,.
$$
This results in an incomplete inner product space; see
\cite{landkof:1972}. Observe also that $L_\infty[0,1]\subset H_\theta 
\subset L_2[0,1]$, and that 
$$
(\psi_1,\psi_2)_\Sigma = (\psi_1,T_\theta\psi_2)\,,
$$
where 
$$
(\psi_1,\psi_2) = \int_0^1 \psi_1(t)\cdot \psi_2(t)\, dt
$$
is the inner product in $L_2[0,1]$, and $T_\theta:\, H_\theta\to
H_\theta$ is defined by 
\begin{equation} \label{e:T}
T_\theta\psi(t) = \int_0^1 \frac{\Sigma\psi(s)}{|t-s|^\theta}ds\,.
\end{equation}

\begin{lemma} \label{l:gaussian.rate.long}
For $\varphi\in L_1[0,1]$ and $1/2<\alpha<1$, 
\begin{equation} \label{e:gauss.rate}
(G_\Sigma)_\alpha^\ast(\varphi) = \sup_{\psi\in L_\infty[0,1]}
  (\psi,\varphi) -
  \frac{\sigma^2}{2}\Bigl(\psi,T_{2\alpha-1}\psi\Bigr)\,, 
\end{equation}
where 
$$
\sigma^2 = (1-\alpha)^2\int_{-\infty}^\infty
|x+1|^{-\alpha}|x|^{-\alpha}
\Bigl[pI_{[x+1\ge 0]}+qI_{[x+1<0]}\Bigr]
\Bigl[pI_{[x\ge 0]}+qI_{[x<0]}\Bigr]\, dx\,,
$$
 $\psi$ is regarded as an element of the dual
space $L_1[0,1]^\prime$, and $T_{2\alpha-1}$ in \eqref{e:T} is
regarded as a map $L_\infty[0,1]\to  L_1[0,1]$. 

(i) Suppose that $\varphi\in T_{2\alpha-1} H_{2\alpha-1}$. Then 
$$
(G_\Sigma)_\alpha^\ast(\varphi) = \frac{1}{2\sigma^2} \| h\|^2_\Sigma\,,
$$
where $\varphi = T_{2\alpha-1}h$. 

(ii) Suppose that Leb$\{ t\in [0,1]:\, \varphi(t)\in K_\Sigma\}>0$,
where $K_\Sigma=Ker(\Sigma)-\{0\}$ is as defined in (\ref{thmsm}). Then
$(G_\Sigma)_\alpha^\ast(\varphi) = \infty$. 
\end{lemma}
\begin{proof}
Note that for $\varphi\in L_1[0,1]$
$$
\int_{-\infty}^\infty 
G_\Sigma\left( \int_0^1 \psi(t)(1-\alpha)|x+t|^{-\alpha}\Bigl[ 
pI_{[x+t\ge 0]}+qI_{[x+t<0]}\Bigr]\, dt\right)
$$
$$
= \frac12(1-\alpha)^2 \int_0^1\int_0^1 \psi(s)\cdot \Sigma\psi(t) \biggl(
\int_{-\infty}^\infty |x+s|^{-\alpha}|x+t|^{-\alpha}
\Bigl[ pI_{[x+s\ge 0]}+qI_{[x+s<0]}\Bigr]
$$
$$
\Bigl[ pI_{[x+t\ge 0]}+qI_{[x+t<0]}\Bigr]\, dx\biggr)ds\,dt
= \frac{\sigma^2}{2} \int_0^1\int_0^1 \frac{\psi(s)\cdot\Sigma
 \psi(t)}{|t-s|^\theta} ds\, dt\,,
$$
and so \eqref{e:gauss.rate} follows. 

For part (i), suppose that $\varphi= T_{2\alpha-1}h$ for $h\in
H_{2\alpha-1}$. For $\psi\in H_{2\alpha-1}$ we have
$$
(\psi,\varphi) - \frac{\sigma^2}{2}(\psi,T_{2\alpha-1}\psi)
= \frac{1}{2\sigma^2}\bigl( h,T_{2\alpha-1}h\bigr) - \frac{\sigma^2}{2}
\Bigl(
(\psi-\frac{1}{\sigma^2}h),T_{2\alpha-1}(\psi-\frac{1}{\sigma^2}h)\Bigr)
$$
because the operator $T_{2\alpha-1}$ is self-adjoint. Therefore,
$$
\sup_{\psi\in  H_{2\alpha-1}} (\psi,\varphi) -
\frac{\sigma^2}{2}(\psi,T_{2\alpha-1}\psi) 
= \frac{1}{2\sigma^2}\bigl( h,T_{2\alpha-1}h\bigr)\,,
$$
achieved at $\psi_0 = h/\sigma^2$, and so by \eqref{e:gauss.rate}, 
$$
(G_\Sigma)_\alpha^\ast(\varphi) \leq \frac{1}{2\sigma^2}\bigl(
h,T_{2\alpha-1}h\bigr)\,. 
$$
On the other hand, for $M>0$ let $\psi_0^{(M)}= \psi_0\one\bigl(
|\psi_0|\leq M\bigr)\in L_\infty[0,1]$. Then
$$
(G_\Sigma)_\alpha^\ast(\varphi) \geq \limsup_{M\to\infty} 
 \psi_0^{(M)}(\varphi) -
  \frac{\sigma^2}{2}\psi_0^{(M)}\Bigl(T_{2\alpha-1}\psi_0^{(M)}\Bigr)
$$
$$ 
= (\psi_0,\varphi) -
\frac{\sigma^2}{2}(\psi_0,T_{2\alpha-1}\psi_0) 
= \frac{1}{2\sigma^2}\bigl( h,T_{2\alpha-1}h\bigr)\,, 
$$
completing the proof of part (i).

For part (ii), note that using \eqref{e:gauss.rate} and choosing for
$c>0$, $\psi(t) = c\varphi(t)/|\varphi(t)|$ if $\varphi(t)\in
K_\Sigma$, and $\psi(t)=0$ otherwise, we obtain 
$$
(G_\Sigma)_\alpha^\ast(\varphi) \geq c\int_A |\varphi(t)|\, dt\,,
$$
where $A= \{ t\in [0,1]:\, \varphi(t)\in K_\Sigma\}$. The proof is
completed by letting $c\to\infty$. 
\end{proof}

\section{Lemmas and  their proofs} \label{sec:lemmas} 
In this section we prove the lemmas used in section \ref{sec:FLDP}. We
retain the notation of  section \ref{sec:FLDP}.

\begin{lemma}\label{l:expon.equiv}
Under any of the assumptions S2, S3, S4, R2, R3 or R4, the families
$\{\mu_n\}$ and $\{\tilde \mu_n\}$ are exponentially 
equivalent in $\mathcal{D}_S$, where $\mathcal{D}$ is the space of all
right-continuous functions with left limits and, as before, the
subscript denotes the sup-norm topology on that space.
\end{lemma}
\begin{proof}
It is clearly enough to consider the case $d=1$. For any $\delta>0$
and $\lambda\in \mathcal{F}_\Lambda\cap -\mathcal{F}_\Lambda$,
$\lambda\not=0$,   
\begin{eqnarray*}
&& \limsup_{n\rightarrow \infty} \frac{1}{b_n} 
\log P\big( ||Y_n-\tilde Y_n||>\delta\big) \\
& \le & 
\limsup_{n\rightarrow \infty} \frac{1}{b_n}\log P\Big(
\frac{1}{a_n}\max_{1\le i\le n}|X_i|>\delta\Big)\\ 
& \le & \limsup_{n\rightarrow \infty} \frac{1}{b_n}
\log \Big( n\, P(|X_1|>a_n\delta) \Big)\\
& \le & \limsup_{n\rightarrow \infty} \frac{1}{b_n}
\Big( \log n -a_n\lambda \delta + \Lambda(\lambda) + \Lambda(-\lambda))
\Big)\\
& = & \limsup_{n\rightarrow \infty} \frac{1}{b_n}\Big( -a_n\lambda
\delta\Big) \,.
\end{eqnarray*}
Under the assumptions  S3, S4,  R3 or R4 we have $a_n/b_n\to\infty$,
so the above limit is equal to $-\infty$. Under the assumptions S2 and
R2, $a_n=b_n$, but we can let $\lambda\rightarrow \infty$ after taking
the limit in $n$. 
\end{proof}

\begin{lemma}\label{expttns}
Under any of the assumptions S2, S3, S4, R2, R3 or R4, the family
$\{\tilde \mu_n\}$ is exponentially 
tight in $\mathcal{D}_S$, i.e, for every $\pi>0$ there exists a
compact  $K_\pi\subset \mathcal{D}_S$, such that 
\[
\lim_{\pi\rightarrow \infty}\limsup_{n\rightarrow \infty
}\frac{1}{b_n}\log \tilde \mu_n(K_\pi^c)=-\infty.
\]  
\end{lemma}

\begin{proof}
We first prove the lemma assuming that  $d=1$. We use the notation
$w(f,\delta)\defn 
\sup\limits_{s,t\in[0,1],|s-t|<\delta}|f(s)-f(t)|$ for the modulus
of 
continuity of a function $f:[0,1]\rightarrow \mathbb{R}^d$. First we
claim that for any $\epsilon>0,$  
\begin{equation}\label{mcont} 
\lim_{\delta\rightarrow 0}\limsup_{n\rightarrow \infty}
\frac{1}{b_n}\log P\big(w(\tilde Y_n,\delta)>\epsilon\big)=-\infty, 
\end{equation} 
where $\tilde Y_n$ is the polygonal process in \eqref{e:polyg}. 
Let us prove the
lemma assuming that the claim is true.  By (\ref{mcont}) and the
continuity of the paths of $\tilde Y_n$, there is $\delta_k>0$ such
that for all $n\geq 1$ 
\[
P\big( w(\tilde Y_n,\delta_k)\ge k^{-1}  \big)\le e^{-\pi b_n   k},
  \] 
and set $A_k=\{f\in
\mathcal{D}:w(f,\delta_k)<k^{-1},f(0)=0\}.$ Now the set $K_\pi\defn
\overline{\cap_{k\ge 1} A_k}$ is compact in $\mathcal{D}_S$ and by
the union of events bound it follows that  
\[  
\limsup_{n\rightarrow \infty }\frac{1}{b_n}\log P(\tilde Y_n\notin
K_\pi)\le -\pi,
\] 
establishing the exponential tightness. 
Next we prove the claim (\ref{mcont}). Observe that for any
$\epsilon>0$, $\delta>0$ small and $n>2/\delta$ 
\begin{eqnarray*}
& P & \big( w(\tilde Y_n,\delta) > \epsilon \big)
  \le  P\Big( \max_{0\le i<j\le n,j-i\le
 [n\delta]+2}\frac{1}{a_n}\Big|\sum_{k=i}^jX_k\Big|>\epsilon \Big)\\ 
& \le & n\sum_{i=1}^{[2n\delta]} P\Big(
 \frac{b_n}{a_n}\Big|\sum_{k=1}^{i}X_k\Big|>b_n\epsilon\Big)\\ 
& \le & ne^{-b_n\lambda \epsilon} \sum_{i=1}^{[2n\delta]}E\Big[ \exp
 \Big( \frac{\lambda b_n}{a_n}\sum_{k=1}^iX_k\Big) +\exp \Big(
 -\frac{\lambda b_n}{a_n}\sum_{k=1}^iX_k\Big) \Big]\\ 
& = & ne^{-b_n\lambda \epsilon} \sum_{i=1}^{[2n\delta]}\Big(\exp\Big[
 \sum_{j\in \mathbb{Z}}\Lambda\Big( \frac{\lambda b_n}{a_n}\phi_{j,i}\Big)
 \Big]+\exp\Big[ \sum_{j\in \mathbb{Z}}\Lambda\Big(-
 \frac{\lambda b_n}{a_n}\phi_{j,i}\Big) \Big]\Big)\\ 
& \le & \frac{2n^2\delta}{e^{b_n\lambda \epsilon}}\Big(\exp\Big[ \sum_{j\in
 \mathbb{Z}}\Lambda\Big( \frac{|\lambda| b_n}{a_n}|\phi|_{j,[2n\delta]}\Big)
 \Big]+\exp\Big[ \sum_{j\in \mathbb{Z}}\Lambda\Big(-
 \frac{|\lambda| b_n}{a_n}|\phi|_{j,[2n\delta]}\Big) \Big]\Big)  
\end{eqnarray*}
by convexity of $\Lambda$ (we use the notation $|\phi|_{i,n} = 
|\phi_{i+1}|+\cdots+|\phi_{i+n}|$ for $i\in \mathbb{Z}$ and $n\ge 1$).
Therefore by lemmas \ref{limits} and \ref{limitr} we have 
$$
 \lim_{\delta\rightarrow 0}\limsup_{n\rightarrow
    \infty}\frac{1}{b_n} \log P\big( w(\tilde Y_n,\delta)>\epsilon
  \big)  \le  -\lambda \epsilon.
$$
Now, letting $\lambda\rightarrow \infty$ we obtain \eqref{mcont}. 

If $d\ge 1$ then $\{\tilde \mu_n\}$ is exponentially tight since
$\{\tilde \mu_n^k\}$, the law of the $k$th coordinate of $\tilde Y_n$,
is exponentially tight for every $1\le k\le d$.  
\end{proof}

\begin{lemma}\label{expttnsl1}
Under the assumptions $S1$ or $R1$ the family $\{\mu_n\}$ is, for any
$p\in [1,\infty)$, 
exponentially tight in the space of functions in $\cap_{p\in [1,\infty)}
  L_p[0,1]$, equipped with the topology $L$, where $f_n$
converges to $f$ if and only if $f_n$ converges to $f$ both pointwise
and in $L_p[0,1]$ for all $p\in [1,\infty)$.  
\end{lemma}
\begin{proof}
Here $a_n=n$ under the assumption $S1$, $a_n=n\Psi_n$ under the
assumption $R1$, and $b_n=n$ in both cases. As before, it is enough to
consider the case $d=1$. We claim that for any $p\in [1,\infty)$, 
\begin{equation}\label{mcont2}
\lim_{x\downarrow 0} \limsup_{n\rightarrow \infty} \frac{1}{n}\log
P\bigg[ \int_0^{1-x}\Bigl|Y_n(t+x)-Y_n(t)\Bigr|^pdt 
\end{equation}
$$
+ \int_0^x \Bigl|Y_n(t)\Bigr|^pdt + 
\int_{1-x}^1 \Bigl|Y_n(t)\Bigr|^pdt >\epsilon\bigg]=-\infty, 
$$
for any $\epsilon>0$, while 
\begin{equation}\label{e:mcont2a}
\lim_{M\uparrow \infty} \limsup_{n\rightarrow \infty} \frac{1}{n}\log
P\Bigl( \sup_{0\leq t\leq 1}|Y_n(t)|>M\Bigr)=-\infty\,.
\end{equation}

Assuming that both claims are true, for any $\pi>0$, $m\geq 1$ and
$k\ge 1,$ we can   choose (using the fact that $Y_n\in  L^\infty[0,1]$
a.s. for all $n\geq 1$)  $0<x_k^{(m)}<1$ such that for all $n\geq 1$,  
\[
P\bigg[ \int_0^{1-x_k^{(m)}}\Bigl|Y_n(t+x_k^{(m)})-Y_n(t)\Bigr|^mdt
\]
\[
+
  \int_0^{x_k^{(m)}}   \Bigl|Y_n(t)\Bigr|^mdt +  
\int_{1-x_k^{(m)}}^1 \Bigl|Y_n(t)\Bigr|^mdt >k^{-1} \Big]\le e^{-\pi
  knm}, 
\]
and $M_\pi>0$ such that for all $n\geq 1$ 
\[
P\Bigl( \sup_{0\leq t\leq 1}|Y_n(t)|>M_\pi\Bigr)\le e^{-\pi n}.
\]
Now define sets  
\[
A_{k,m}= \Big\{f\in \cap_{p\geq   1}L_p[0,1]:\int_0^{1-x_k^{(m)}}
\Bigl|f(t+x_k^{(m)})-f(t)\Bigr|^mdt 
\] 
$$
+ \int_0^{x_k^{(m)}}
\Bigl|f(t)\Bigr|^mdt +  \int_{1-x_k^{(m)}}^1 \Bigl|f(t)\Bigr|^mdt
\le k^{-1},\ \sup_{0\leq t\leq 1}|f(t)|\le
M_\pi\Big\},
$$
and set $K_\pi=\overline{\cap_{k,m\ge 1}A_{k,m}}.$ Then $K_\pi$ is
compact for every $\pi>0$ by Tychonov's theorem (see 
 theorem 19, p. 166 in \cite{royden:1968} and theorem 20, p. 298 in
 \cite{dunford:schwartz:1988}). Furthermore, 
\[
\limsup_{n\rightarrow \infty} \frac{1}{n}\log P[Y_n\notin K_\pi]\le
-\pi.
\]

This will complete the proof once we prove (\ref{mcont2}) and
\eqref{e:mcont2a}. We first prove (\ref{mcont2}) for $p=1$. 
Observe that 
\begin{eqnarray*}
&&P\Big[ \int_0^{1-x}|Y_n(t+x)-Y_n(t)|dt>\epsilon \Big] \leq 
 P\Big[\frac{[nx]}{n} 
    \frac{1}{a_n}\sum_{i=1}^n|X_i|>\epsilon\Big]\\
&\le &  e^{-\lambda n \epsilon/x}E\Big[ \exp \Big( \lambda
    \frac{b_n}{a_n}\sum_{i=1}^n|X_i|\Big)\Big] 
 \le  e^{-\lambda n \epsilon/x} E\Big[ \prod_{i=1}^n
    \exp\Big(\frac{\lambda b_n}{a_n} |X_i| \Big) \Big]\\
& \le&  e^{-\lambda n \epsilon/x}  E\Big[ \prod_{i=1}^n \Big(
    \exp\Big(\frac{\lambda b_n}{a_n} X_i \Big)
    +\exp\Big(-\frac{\lambda b_n}{a_n}X_i \Big)\Big)\Big]\\
& = & e^{-\lambda n \epsilon/x}\sum_{l_i=\pm 1}E\Big[ \exp\Big(
    \frac{\lambda b_n}{a_n} \sum_{i=1}^n l_iX_i\Big)\Big]\\
& = & e^{-\lambda n \epsilon/x} \sum_{l_i=\pm 1}\exp\Big( \sum_{j\in
    \mathbb{Z}}
  \Lambda\Big(\frac{\lambda b_n}{a_n}(\phi_{j+1}l_1+\cdots+\phi_{j+n}l_n)
  \Big) \Big)\\
&\le&  2^n e^{-\lambda n \epsilon/x} \exp\Big( \sum_{j\in
    \mathbb{Z}}\Lambda\Big(\frac{\lambda b_n}{a_n}|\phi|_{j,n}
  \Big)+\sum_{j\in\mathbb{Z}}\Lambda\Big(
  -\frac{\lambda b_n}{a_n}|\phi|_{j,n}\Big) \Big) \,.
\end{eqnarray*}
Therefore, 
\begin{eqnarray*}
&&\limsup_{n\rightarrow \infty} \frac{1}{n} \log P\Big[
    \int_0^{1-x}|Y_n(t+x)-Y_n(t)|dt>\epsilon \Big] \\  
& \le & \log 2 -\frac{\lambda\epsilon}{x}  +\limsup_{n\rightarrow
    \infty }\frac{1}{n} \sum_{j\in
    \mathbb{Z}}\Lambda\Big(\frac{\lambda b_n}{a_n}|\phi|_{j,n}
    \Big)+\limsup_{n\rightarrow \infty }\frac{1}{n} \sum_{j\in
    \mathbb{Z}}\Lambda\Big(-\frac{\lambda b_n}{a_n}|\phi|_{j,n} \Big)\,.
\end{eqnarray*}
Keeping $\lambda>0$ small, using  lemma \ref{limits} and lemma
\ref{limitr} and then letting $x\rightarrow 0$ one establishes the limit 
\[ 
 	\limsup_{n\rightarrow \infty} \frac{1}{n} \log P\Big[
    \int_0^{1-x}|Y_n(t+x)-Y_n(t)|dt>\epsilon \Big] =- \infty.   
\]
It is simpler to show a similar inequality for the 
second and the third integrals under the probability of the equation \eqref{mcont2}. The proof of
\eqref{e:mcont2a} is similar, starting with
$$
P\Bigl( \sup_{0\leq t\leq 1}|Y_n(t)|>M\Bigr)
\leq P\Big( \frac{1}{a_n}\sum_{i=1}^n|X_i|>M\Big)\,.
$$
Now one establishes (\ref{mcont2}) for $p\ge 1$ by writing, for $M>0$,
\pagebreak
$$
P\bigg[ \int_0^{1-x}\Bigl|Y_n(t+x)-Y_n(t)\Bigr|^pdt 
+ \int_0^x \Bigl|Y_n(t)\Bigr|^pdt + 
\int_{1-x}^1 \Bigl|Y_n(t)\Bigr|^pdt >\epsilon\bigg]
$$
$$
\leq P\bigg[ \int_0^{1-x}\Bigl|Y_n(t+x)-Y_n(t)\Bigr|dt 
+ \int_0^x \Bigl|Y_n(t)\Bigr|dt + 
\int_{1-x}^1 \Bigl|Y_n(t)\Bigr|dt >\frac{\epsilon}{2M^{p-1}}\bigg]
$$
$$ 
+ P\Bigl[ \sup_{0\leq t\leq 1}|Y_n(t)|>M\Bigr]\,,
$$
and letting first $n\to\infty$, $x\downarrow 0$, and then
$M\uparrow\infty$. 
\end{proof}

\begin{lemma}\label{l:spprt.BV}
Under the assumptions $S1$ or $R1$, the corresponding upper rate
functions, $G^{sl}$ in \eqref{upperratefn} and $G^{rl}$ 
in \eqref{e:upperratefn.1}, are infinite outisde of
the space $\mathcal{BV}$. 
\end{lemma}
\begin{proof}
Let $f\notin\mathcal{BV}$. Choose $\delta>0$ small enough such that any $\lambda$
with $|\lambda|\leq \delta$ is in $\mathcal{F}_\Lambda^\circ$ and a
vector with $k$ identical components $(\lambda,\ldots ,\lambda)$ is in
the interiors of both $\Pi_{t_1,\ldots,t_k}$ in \eqref{pi} and
$\Pi^{r,\alpha}_{t_1,\ldots,t_k}$ in \eqref{pir} and \eqref{pir.1}. 
For $M>0$ choose a partition
$0<t_1<\cdots<t_k= 1$ of $[0,1]$ such that
$\sum_{i=1}^{k}\bigl|f(t_i) -f(t_{i-1})\bigr|>M$. For $i=1,\ldots, k$
such that $f(t_i) -f(t_{i-1})\not=0$ choose $\lambda_i$ of length
$\delta$ in the direction of $f(t_i) -f(t_{i-1})$. Then under, say,
assumption $S1$, 
$$
G^{sl}(f)\geq \sup_{\underline{\lambda}\in \Pi_{t_1,\ldots,t_k} }
\sum_{i=1}^{k}\Big\{\lambda_i\cdot\big(f(t_i)
-f(t_{i-1})\big)-(t_i-t_{i-1})\Lambda(\lambda_i)\Big\}
$$
$$
\geq \delta M - \sup_{|\lambda|\leq \delta} \Lambda(\lambda)\,.
$$
Letting $M\to\infty$ proves the statement under the assumption $S1$,
and the argument under the assumption $R1$ is similar. 
\end{proof}

\begin{lemma}\label{limits}
Suppose $\Lambda:\mathbb{R}^d\rightarrow \mathbb{R}$ is  the
log-moment generating function of a mean zero random variable $Z$,  with
$0\in \mathcal{F}_\Lambda^\circ$,
$\sum\limits_{i=-\infty}^\infty|\phi_i|<\infty$ with
$\sum\limits_{i=-\infty}^\infty\phi_i=1$ and    $0< t_1<\cdots< t_k\le
1$. 
\begin{enumerate}[(i)]
\item  For all $\underline{\lambda}=(\lambda_1,\ldots,\lambda_k)\in
  \Pi_{t_1,\ldots,t_k}\subset (\mathbb{R}^{d})^k$, 
\[\lim_{n\rightarrow \infty} \frac{1}{n}\sum_{l=-\infty}^\infty
  \Lambda\Big(
  \sum_{i=1}^{k}\lambda_i\phi_{l+[nt_{i-1}],[nt_i]-[nt_{i-1}]}
  \Big)=\sum_{i=1}^k(t_i-t_{i-1})\Lambda(\lambda_i).\] 

\item If $a_n/\sqrt{n}\rightarrow \infty$ and $a_n/n\rightarrow 0$
  then for all  $\underline{\lambda}\in (\mathbb{R}^d)^k$,  
\[ \lim_{n\rightarrow \infty} \frac{n}{a_n^2}\sum_{l=-\infty}^\infty
  \Lambda\Big(
  \frac{a_n}{n}\sum_{i=1}^{k}\lambda_i\phi_{l+[nt_{i-1}],[nt_i]-[nt_{i-1}]}
  \Big)=\sum_{i=1}^k(t_i-t_{i-1})\lambda_i\cdot \Sigma\lambda_i,\] 
where $\Sigma$ is the covaraince matrix of $Z$. 

\item If  $\Lambda(\cdot)$ is balanced regular varying at $\infty$
  with exponent $\beta>1$, $a_n/n\rightarrow \infty$ and $b_n$ is as
  defined  as defined in assumption $S4$, then for all
  $\underline{\lambda}\in (\mathbb{R}^d)^k$,  
\[ 
\lim_{n\rightarrow \infty} \frac{1}{b_n}\sum_{l=-\infty}^\infty
  \Lambda\Big(
  \frac{b_n}{a_n}\sum_{i=1}^{k}\lambda_i\phi_{l+[nt_{i-1}],[nt_i]-[nt_{i-1}]}
  \Big)
\]
\[
=\sum_{i=1}^k(t_i-t_{i-1})\zeta\Big(\frac{\lambda_i}{|\lambda_i|}
  \Big)|\lambda_i|^\beta. \] 
\end{enumerate}
\end{lemma}

\begin{proof}
(i) We begin  by making a few observations:

\begin{enumerate}[$(a)$]
\item For every $\delta>0$ there exists $N_\delta$ such that for all
  $n>N_\delta$ 
\begin{equation}\label{i}
  \sum_{|i|>(n\min\limits_j(t_j-t_{j-1}))^{1/2}}|\phi_i|<\delta.
\end{equation} 

\item  For fixed
  $\underline{\lambda}=(\lambda_1,\ldots,\lambda_k)\in
  \Pi_{t_1,\ldots,t_k}$, there exists $M>0$ such that  for all
  $l\in\mathbb{Z}$ and all $n$ large enough 
\begin{equation}\label{ii}
\Big|\Lambda\Big(\sum_{i=1}^{k}\lambda_i\phi_{l+[nt_{i-1}],s_i}\Big)
\Big|\le M,
\end{equation}where
$s_i=s_i(n)=[nt_i]-[nt_{i-1}]$.  Since the zero mean of $Z$ means that $\Lambda(x) = o(|x|)$ as $|x|\to 0$, it follows from \eqref{ii} that 
there exists $C>0$ such that in the same range of $n$ and for all
$l\in\mathbb{Z}$ 
\begin{equation}\label{iii}
\Big|\Lambda\Big(\sum_{i=1}^{k}\lambda_i\phi_{l+[nt_{i-1}],s_i}\Big)
\Big|\le C\Big| \sum_{i=1}^{k}\lambda_i\phi_{l+[nt_{i-1}],s_i}
\Big|\,.
\end{equation}
\end{enumerate}

\noindent Let $L=\Big(|\lambda_1|+\cdots+|\lambda_k|\Big)$. Since
$\Lambda$ is continuous at $\lambda_j$, given $\epsilon>0$ we can
choose $\delta>0$ so that for $n$ large enough, 
\[
\Big|\sum_{i=1}^{k}\lambda_i\phi_{l+[nt_{i-1}],s_i}-\lambda_j\Big|<\delta 
\]
for all $-[nt_j]+\sqrt{s_j}<l<-[nt_{j-1}]-\sqrt{s_j}$, and then 
\[
\Big|\frac{1}{n} \sum_{l=-[nt_j]+\sqrt{s_j}}^{-[nt_{j-1}]-\sqrt{s_j}}
\Lambda\Big(\sum_{i=1}^{k}\lambda_i\phi_{l+[nt_{i-1}],s_i}\Big)
-\frac{s_j-2\sqrt{s_j}}{n} \Lambda(\lambda_j) \Big|<\epsilon. 
\]
Therefore for $j=1,\ldots, k$ 
\begin{equation}\label{part1}
\lim_{n\rightarrow\infty}\frac{1}{n}
\sum_{l=-[nt_j]+\sqrt{s_j}}^{-[nt_{j-1}]-\sqrt{s_j}}
\Lambda\Big(\sum_{i=1}^{k}\lambda_i\phi_{l+[nt_{i-1}],s_i}\Big)=
(t_j-t_{j-1})\Lambda(\lambda_j). 
\end{equation}
Note that 
\begin{equation}\label{part2}
\Big|\frac{1}{n} \sum_{l=-[nt_j]-\sqrt{s_j}}^{-[nt_j]+\sqrt{s_{j+1}}} 
\Lambda\Big(\sum_{i=1}^{k} \lambda_i\phi_{ l+[nt_{i-1}] ,s_i }  \Big)
\Big|\stackrel{(\ref{ii})}{\le} \frac{\sqrt{s_j}+\sqrt{s_{j+1}}}{n}M 
\stackrel{n\rightarrow \infty}{\longrightarrow}  0\,.
\end{equation}
Finally, observe that for large $n$,
\begin{eqnarray}\label{part3}
&& \Big|\frac{1}{n} \sum_{l=-\infty}^{-[nt_k]-\sqrt{s_k}} \Lambda
\Big(\sum_{i=1}^{k}\lambda_i\phi_{l+[nt_{i-1}],s_i} \Big)
\Big|\nonumber\\
&\stackrel{(\ref{iii})}{\le} & C\frac{1}{n} 
\sum_{l=-\infty}^{-[nt_k]-\sqrt{s_k}} \Big|\sum_{i=1}^{k}\lambda_i
\phi_{l+[nt_{i-1}],s_i} \Big|\nonumber\\
& \le & CL\sum_{l=-\infty}^{-\sqrt{s_k}}|\phi_l|
\stackrel{(i)}{\rightarrow} 0.
\end{eqnarray}
and 
\begin{eqnarray}\label{part4}
&& \Big|\frac{1}{n} \sum_{l=\sqrt{s_1}}^{\infty}
  \Lambda\Big(\sum_{i=1}^{k} \lambda_i\phi_{l+[nt_{i-1}],s_i} \Big)
\Big|\nonumber\\
&\stackrel{(\ref{iii})}{\le} & C\frac{1}{n} 
  \sum_{l=\sqrt{s_1}}^{\infty}
 \Big|\sum_{i=1}^{k}\lambda_i\phi_{l+[nt_{i-1}],s_i} \Big|\nonumber\\
& \le & CL\sum_{l=\sqrt{s_1}}^{\infty}|\phi_l| \rightarrow 0.
\end{eqnarray}
Thus, combining (\ref{part1}), (\ref{part2}), (\ref{part3}) and
(\ref{part4}) we have 
\[
\lim_{n\rightarrow \infty} \frac{1}{n} \sum_{l=-\infty}^{\infty}
\Lambda \Big(\sum_{i=1}^{k}\lambda_i
\phi_{l+[nt_{i-1}],[nt_{i}]-[nt_{i-1}]}\Big)
= \sum_{i=1}^{k}(t_i-t_{i-1})\Lambda(\lambda_i).
\] 

(ii) Since $\Lambda(x) \sim x\cdot\Sigma x/2$ as 
$|x|\rightarrow 0$, we see that for every $1\le j\le k$, 
\[
\lim_{n\rightarrow\infty}\frac{n}{a_n^2} 
\sum_{l=-[nt_j]+\sqrt{s_j}}^{-[nt_{j-1}]-\sqrt{s_j}} 
\Lambda\Big(\frac{a_n}{n}\sum_{i=1}^{k}\lambda_i
\phi_{l+[nt_{i-1}], [nt_{i}]-[nt_{i-1}]}\Big)=(t_j-t_{j-1})
 \frac{1}{2}\lambda_j\cdot \Sigma \lambda_j.
\] 
The rest of the proof is similar to the proof of part (i).

(iii)  Since $\Lambda(\lambda)$ is regular varying at infinity with
exponent $\beta>1$, for every $1\le j\le k$, 
\pagebreak
\[
\lim_{n\rightarrow\infty}\frac{1}{b_n}
\sum_{l=-[nt_j]+\sqrt{s_j}}^{-[nt_{j-1}]-\sqrt{s_j}}
\Lambda\Big(\frac{b_n}{a_n}\sum_{i=1}^{k}\lambda_i\phi_{l+[nt_{i-1}],
  [nt_{i}]-[nt_{i-1}]}\Big) 
\]
\[
=(t_j-t_{j-1})\zeta
\Big(\frac{\lambda_j}{|\lambda_j|} \Big)|\lambda_j|^\beta.
\] 
The rest of the proof is, once again, similar to the proof of part (i). 
\end{proof}

\begin{lemma}\label{limitr}
Suppose $\Lambda:\mathbb{R}^d\rightarrow \mathbb{R}$ is  the
log-moment generating function of a mean zero random variable,  with
$0\in \mathcal{F}_\Lambda^\circ$, the coefficients of the
moving average are balanced regularly varying with exponent $\alpha$
as in Assumption \ref{regularvarying}, and $0< t_1<\cdots< t_k\le 1$. 
\begin{enumerate}[(i)]
\item  For all $\underline{\lambda}=(\lambda_1,\ldots,\lambda_k)\in
  \Pi^{r,\alpha}_{t_1,\ldots,t_k}\subset (\mathbb{R}^{d})^k$, 
\[
\lim_{n\rightarrow \infty} \frac{1}{n}\sum_{l=-\infty}^\infty
\Lambda\Big(
\frac{1}{\Psi_n}\sum_{i=1}^{k}\lambda_i\phi_{l+[nt_{i-1}],[nt_i]-[nt_{i-1}]}
\Big)=\Lambda_{t_1,\cdots,t_k}^{rl}(\underline{\lambda}).
\] 

\item If $a_n/\sqrt{n}\rightarrow \infty$ and $a_n/n\rightarrow 0$
  then  for all   $\underline{\lambda}\in (\mathbb{R}^d)^k$,  
$$
\lim_{n\rightarrow \infty}
\frac{n\Psi_n^2}{a_n^2}\sum_{l=-\infty}^\infty 
\Lambda\Big( \frac{a_n}{n\Psi_n^2}\sum_{i=1}^{k}\lambda_i
\phi_{l+[nt_{i-1}],[nt_i]-[nt_{i-1}]}
\Big)
$$
$$
=\left\{ \begin{array}{lcl}\int\limits_{-\infty}^\infty
   G_\Sigma\Big(h_{t_1,\ldots,t_k}(x;\underline{\lambda}) \Big) dx & if
   & \alpha<1\\ \sum\limits_{i=1}^{k}(t_i-t_{i-1})G_\Sigma(\lambda_i) &
   if & \alpha=1,\end{array}\right..
$$

\item If $a_n/n\rightarrow \infty$, $b_n$ is as defined in assumption
$R4$,  and $\Lambda(\cdot)$ is balanced regular varying at $\infty$
with exponent $\beta>1$, then for all
$\underline{\lambda}\in (\mathbb{R}^d)^k$,  
$$
\lim_{n\rightarrow \infty} \frac{1}{b_n}\sum_{l=-\infty}^\infty
\Lambda\Big(
\frac{b_n}{a_n}\sum_{i=1}^{k}\lambda_i\phi_{l+[nt_{i-1}],[nt_i]-[nt_{i-1}]}
\Big) 
$$
$$
=\left\{ \begin{array}{lcl}\int\limits_{-\infty}^\infty
   \Lambda^h\Big(h_{t_1,\ldots,t_k}(x;\underline{\lambda}) \Big) dx & if
   & \alpha<1\\ \sum\limits_{i=1}^{k}(t_i-t_{i-1})\Lambda^h(\lambda_i) &
   if & \alpha=1,\end{array}\right..
$$
\end{enumerate}
\end{lemma}

\begin{proof}
(i) We may (and will) assume that $t_k=1$, since we can always add an
  extra point with the zero vector $\lambda$ corresponding to it. 
Let us first assume that $\alpha<1$. Note that for any $m\ge 1$
  and large $n$, 
\begin{eqnarray*}
&&\frac{1}{n}\sum_{j=nm+1}^{n(m+1)} \Lambda\Big( \frac{1}{\Psi_n}
  \sum_{i=1}^{k}\lambda_i\phi_{j+[nt_{i-1}],[nt_i]-[nt_{i-1}]}\Big)\\ 
& = & \frac{1}{n}\sum_{j=nm+1}^{n(m+1)}
  \Lambda\Big(\sum_{i=1}^k\lambda_i\frac{n\psi(n)}{\Psi_n}
  \frac{1}{n}\big( \frac{\phi_{j+[nt_{i-1}]+1}}{\psi(n)}+\cdots+
\frac{\phi_{j+[nt_i]}}{\psi(n)} \big)\Big) \\
& = & \int_m^{m+1} f_n(x)\, dx\,,
\end{eqnarray*}
where 
$$
f_n(x) = \Lambda\Big( \frac{1}{\Psi_n}
  \sum_{i=1}^{k}\lambda_i\phi_{j+[nt_{i-1}],[nt_i]-[nt_{i-1}]}\Big)
$$
if $(j-1)/n<x\leq j/n$ for $j=nm+1,\ldots ,n(m+1)$. 

Notice that by Karamata's theorem (see \cite{resnick:1987}),
$n\psi(n)/\Psi_n\to 1-\alpha$ as $n\to\infty$. Furthermore, given
$0<\epsilon<\alpha$, we can use Potter's bounds (see Proposition 0.8 {\it
  ibid}) to check that there is $n_\epsilon$ such that for all $n\geq
n_\epsilon$, for all $k=[nt_{i-1}]+1, \ldots, [nt_i]$, $m-1<x\leq m$
and $(j-1)/n<x\leq j/n$ 
$$
\frac{\phi_{j+k}}{\psi(n)} =
\frac{\phi_{j+k}}{\psi(j+k)}\frac{\psi(j+k)}{\psi(j)}
\frac{\psi(j)}{\psi(n)} 
$$
$$
\in \left( (1-\epsilon)\, p\left(
\frac{j+k}{j}\right)^{-(\alpha+\epsilon)}x^{-\alpha}, \, 
(1+\epsilon)\, p\left(
\frac{j+k}{j}\right)^{-(\alpha-\epsilon)}x^{-\alpha}\right),
$$
and so for $n$ large enough,
\begin{equation} \label{e:potter}
\frac{1}{n}\left( \frac{\phi_{j+[nt_{i-1}]+1}}{\psi(n)}+\cdots+
\frac{\phi_{j+[nt_i]}}{\psi(n)} \right)
\end{equation}
$$
\in \left( (1-\epsilon)\, p \int_{t_{i-1}}^{t_i} \left(
\frac{y+x}{x}\right)^{-(\alpha+\epsilon)}x^{-\alpha}\, dy, \,
(1+\epsilon)\, p \int_{t_{i-1}}^{t_i} \left(
\frac{y+x}{x}\right)^{-(\alpha-\epsilon)}x^{-\alpha}\, dy \right)\,.
$$
Therefore,
$$
\frac{1}{\Psi_n}
  \sum_{i=1}^{k}\lambda_i\phi_{j+[nt_{i-1}],[nt_i]-[nt_{i-1}]}
\to (1-\alpha)\, p \sum_{i=1}^k \lambda _i\int_{t_{i-1}}^{t_i}
(y+x)^{-\alpha}\, dy 
$$
$$
= p \sum_{i=1}^k \lambda _i \Bigl(
(t_i+x)^{1-\alpha} - (t_{i-1}+x)^{1-\alpha}\Bigr). 
$$
This last vector is a convex linear combination of the vectors 
$p\bigl( (1+x)^{1-\alpha} - x^{1-\alpha}\bigr)\lambda_i$, $i=1\ldots,
k$. By the definition of the set $\Pi^{r,\alpha}_{t_1,\ldots,t_k}$,
each one of these vectors belongs to $\mathcal{F}_\Lambda^\circ$ and,
by convexity of $\Lambda$, so does the convex linear
combination. Therefore, 
$$
\Lambda\Big( \frac{1}{\Psi_n}
  \sum_{i=1}^{k}\lambda_i\phi_{j+[nt_{i-1}],[nt_i]-[nt_{i-1}]}\Big)
\to \Lambda\Bigl( p \sum_{i=1}^k \lambda _i \Bigl(
(t_i+x)^{1-\alpha} - (t_{i-1}+x)^{1-\alpha}\Bigr)\Bigr). 
$$
This convexity argument also shows that the function $f_n$ is
uniformly bounded on $(m, m+1]$ for large enough $n$, and so we
conclude that for any $m\ge 1$
$$
\frac{1}{n}\sum_{j=nm+1}^{n(m+1)} \Lambda\Big( \frac{1}{\Psi_n}
  \sum_{i=1}^{k}\lambda_i\phi_{j+[nt_{i-1}],[nt_i]-[nt_{i-1}]}\Big)
$$
$$
\to  \int_m^{m+1} \Lambda\Big(
  (1-\alpha)\sum_{i=1}^k\lambda_i\int\limits_{x+t_{i-1}}^{x+t_{i}}p
  y^{-\alpha}dy\Big) dx. 
$$
Similar arguments show that for $m\le -3$
$$
\frac{1}{n}\sum_{j=nm+1}^{n(m+1)} \Lambda\Big(
  \frac{1}{\Psi_{n}}\sum_{i=1}^{k}\lambda_i
  \phi_{j+[nt_{i-1}],[nt_i]-[nt_{i-1}]}\Big)
$$
$$
\to  \int_m^{m+1} \Lambda\Big(
  (1-\alpha)\sum_{i=1}^k\lambda_i\int\limits_{x+t_{i-1}}^{x+t_{i}}q
  |y|^{-\alpha}dy\Big) dx,   
$$
and that for any $\delta>0$, 
\[ 
\frac{1}{n}\sum_{j=-2n+1}^{-n-n\delta} \Lambda\Big(
 \frac{1}{\Psi_{n}}\sum_{i=1}^{k}\lambda_i
 \phi_{j+[nt_{i-1}],[nt_i]-[nt_{i-1}]}\Big) 
\]
\[
\rightarrow
 \int_{-2}^{-1-\delta} \Lambda\Big(
 (1-\alpha)\sum_{i=1}^k\lambda_i\int\limits_{x+t_{i-1}}^{x+t_{i}}q
 |y|^{-\alpha}dy\Big) dx
\] 
 and 
\[ 
\frac{1}{n}\sum_{j=n\delta}^{n} \Lambda\Big(
\frac{1}{\Psi_{n}}\sum_{i=1}^{k}\lambda_i
\phi_{j+[nt_{i-1}],[nt_i]-[nt_{i-1}]}\Big) 
\]
\[
\rightarrow
\int_{\delta}^{1} \Lambda\Big(
(1-\alpha)\sum_{i=1}^k\lambda_i\int\limits_{x+t_{i-1}}^{x+t_{i}}
py^{-\alpha}dy\Big) dx.
\]
Using once again the same argument we see that for small $\delta$ 
\pagebreak
$$
\frac{1}{n}\sum^{0}_{j=-n} \one\Bigl(
  \Bigl|\frac{j}{n}+t_i\Bigr|>\delta \ \text{all} \ i=1,\ldots, k\Bigr)
\Lambda\Big( \frac{1}{\Psi_n}
  \sum_{i=1}^{k}\lambda_i\phi_{j+[nt_{i-1}],[nt_i]-[nt_{i-1}]}\Big)
$$
$$
 \to  \int_{-1}^{0} \one\Bigl(
  |x+t_i|>\delta \ \text{all} \ i=1,\ldots, k\Bigr)
$$
$$
\Lambda\Big( (1-\alpha)
\sum_{i=1}^k\lambda_i\int\limits_{x+t_{i-1}}^{x+t_{i}}|y|^{-\alpha}\big(pI_{[y\ge
  0]}+qI_{[y<0]} \big)dy\Big) dx. 
$$
We have covered above all choices of the subscript $j$ apart from a
finite number of stretches of $j$ of length at most $n\delta$ each. By
the definition of the set $\Pi^{r,\alpha}_{t_1,\ldots,t_k}$ we see
that there is a finite $K$ such that for all $n$ large enough,  
\[
\frac{1}{n}\sum_{j \, \text{not yet considered}}  \Lambda\Big(
\frac{1}{\Psi_{n}}\sum_{i=1}^{k}\lambda_i
\phi_{j+[nt_{i-1}],[nt_i]-[nt_{i-1}]}\Big)\le K\delta.
\] 
It follows from \eqref{e:potter} and the fact that
$\Lambda(\lambda)=O(|\lambda|^2)$ as  $\lambda\rightarrow 0$ that for
all $|m|$ large enough there is $C\in (0,\infty)$ such that
$$
\frac{1}{n}\sum^{n(m+1)}_{nm+1} \Lambda\Big( \frac{1}{\Psi_n}
\sum_{i=1}^{k}\lambda_i\phi_{j+[nt_{i-1}],[nt_i]-[nt_{i-1}]}\Big) 
\leq C|m|^{-2\alpha}
$$
for all $n$ large enough. This is summable by the assumption on
$\alpha$, and so the dominated convergence theorem gives us the result. 

Next we move our attention to the case $\alpha=1$. Choose any
$\delta>0$. By the slow variation of $\Psi_n$ we see that
\[
\sup_{j>\delta n \, \text{or}\, j<-(1+\delta)n}
\frac{|\phi_{j,n}|}{\Psi_{n}}\rightarrow 0\,,
\] 
while for any $0<x<1$ we
have  
\[ 
\frac{\phi_{0,[nx]}}{\Psi_{n}}\rightarrow p \mbox{ and  }
\frac{\phi_{-[nx],[nx]}}{\Psi_{n}}\rightarrow q.
\] 
Write
\[
\frac{1}{n}\sum_{j=-n+1}^0\Lambda\Big( \frac{1}{\Psi_n}
\sum_{i=1}^{k}\lambda_i\phi_{j+[nt_{i-1}],[nt_i]-[nt_{i-1}]}\Big)
\]
\[
= \sum_{m=1}^k \frac{1}{n} \sum_{j=-[nt_m]+1}^{j=-[nt_{m-1}]}
\Lambda\Big( \sum_{i=1}^{k}\lambda_i
\frac{\phi_{j+[nt_{i-1}],[nt_i]-[nt_{i-1}]}}{\Psi_n}\Big) \,.
\]
Fix $m=1,\ldots, k$, and observe that for any $\epsilon>0$ and $n$
large enough,
$$
\frac{1}{n} \sum_{j=-[nt_m]+1}^{-[nt_{m-1}]}
\Lambda\Big( \sum_{i=1}^{k}\lambda_i
\frac{\phi_{j+[nt_{i-1}],[nt_i]-[nt_{i-1}]}}{\Psi_n}\Big)
= \int\limits_{-t_m-\epsilon}^{-t_{m-1}} f_n(x)\, dx\,,
$$
where this time
$$
f_n(x) = \one\Bigl( -\frac{[nt_m]}{n}<x\leq
-\frac{[nt_{m-1}]}{n}\Bigr)
\Lambda\Big( \sum_{i=1}^{k}\lambda_i
\frac{\phi_{j+[nt_{i-1}],[nt_i]-[nt_{i-1}]}}{\Psi_n}\Big)
$$
if  $(j-1)/n<x\leq j/n$ for $j=-[nt_m]+1,\ldots ,-[nt_{m-1}]$,
otherwise $f_n(x)=0$. Clearly, $f_n(x)\to 0$ as $n\to\infty$ for all
$-t_m-\epsilon< x<-t_m$. Furthermore, 
$$
\frac{\phi_{j+[nt_{i-1}],[nt_i]-[nt_{i-1}]}}{\Psi_n}\to 0
$$
uniformly in $i\not= m$ and $j=-[nt_m]+1,\ldots ,-[nt_{m-1}]$, while 
for every $-t_m<x<-t_{m-1}$, 
$$
\frac{\phi_{j+[nt_{m-1}],[nt_i]-[nt_{m-1}]}}{\Psi_n}\to p+q=1\,.
$$
By the definition of the set $\Pi^{r,1}_{t_1,\ldots,t_k}$ we see that 
$f_n \to \one_{(-t_{m},-t_{m-1})}\Lambda(\lambda_m)$ a.e., and that
the functions $f_n$ are uniformly bounded for large $n$. Therefore, 
$$
\frac{1}{n}\sum_{j=-n+1}^0\Lambda\Big( \frac{1}{\Psi_n}
\sum_{i=1}^{k}\lambda_i\phi_{j+[nt_{i-1}],[nt_i]-[nt_{i-1}]}\Big)
\to \sum_{m=1}^k (t_m-t_{m-1})\Lambda(\lambda_m)\,.
$$
Finally, the argument above, using Potter's bounds and the fact that
$\Lambda(\lambda) = O(|\lambda|^2)$ as $\lambda\to 0$, shows that 
\[
\frac{1}{n}\sum_{j\notin[-n,0]} \Lambda\Big( \frac{1}{\Psi_n}
\sum_{i=1}^{k}\lambda_i\phi_{j+[nt_{i-1}],[nt_i]-[nt_{i-1}]}\Big)
\rightarrow 0. 
\] 
This completes the proof of part (i). 

For part (ii) consider, once again, the cases $1/2<\alpha<1$ and
$\alpha=1$ separately. If $1/2<\alpha<1$, then for every $m\geq 1$ we
use the regular variation and the fact that $\Lambda(x) \sim
x\cdot\Sigma x/2$ as  $|x|\rightarrow 0$ to obtain 
$$
\frac{n\Psi_n^2}{a_n^2}\sum_{j=nm+1}^{n(m+1)}
\Lambda\Big( \frac{a_n}{n\Psi_n^2}\sum_{i=1}^{k}\lambda_i
\phi_{l+[nt_{i-1}],[nt_i]-[nt_{i-1}]} \Big)
\to 
$$
$$
\int\limits_m^{m+1}
\Bigl( (1-\alpha)\sum_{i=1}^k\lambda_i\int\limits_{x+t_{i-1}}^{x+t_{i}}p 
  y^{-\alpha}dy\Bigr) \cdot\Sigma  \Bigl(
  (1-\alpha)\sum_{i=1}^k\lambda_i\int\limits_{x+t_{i-1}}^{x+t_{i}}p
  y^{-\alpha}dy\Bigr)/2\, dx\,,
$$
and we proceed as in the proof of part (i), considering the various
other ranges of $m$, obtaining the result. If $\alpha=1$, then for any 
$m=1,\ldots, k$, by the regular variation and the fact that
$\Lambda(x) \sim x\cdot\Sigma x/2$ as  $|x|\rightarrow 0$, one has 
$$
\frac{n\Psi_n^2}{a_n^2} \sum_{j=-[nt_m]+1}^{[nt_{m-1}]}
\Lambda\Big( \frac{a_n}{n\Psi_n}\sum_{i=1}^{k}\lambda_i
\frac{\phi_{j+[nt_{i-1}],[nt_i]-[nt_{i-1}]}}{\Psi_n}\Big)
\to \int\limits_{-t_m}^{-t_{m-1}} \frac12 \lambda_m\cdot \Sigma \lambda_m\,
dx\,,
$$
and so 
$$
\frac{n\Psi_n^2}{a_n^2} \sum_{j= -n+1}^{0}
\Lambda\Big( \frac{a_n}{n\Psi_n^2}\sum_{i=1}^{k}\lambda_i
\phi_{j+[nt_{i-1}],[nt_i]-[nt_{i-1}]}\Big)
\to \frac12 \sum_{m=1}^k (t_m-t_{m-1})\lambda_m\cdot \Sigma
\lambda_m\,.
$$
As in part (i), by using Potter's bounds and the fact that
$\Lambda(\lambda) = O(|\lambda|^2)$ as $\lambda\to 0$, we see that 
\[
\frac{n\Psi_n^2}{a_n^2}\sum_{j\notin[-n,0]} \Lambda\Big(
\frac{a_n}{n\Psi_n^2} 
\sum_{i=1}^{k}\lambda_i\phi_{j+[nt_{i-1}],[nt_i]-[nt_{i-1}]}\Big)
\rightarrow 0, 
\] 
giving us the desired result. 

We proceed in a similar fashion in part (iii). If $1/2<\alpha<1$, then,
for example, for $m\geq 1$, by the regular variation at infinity, 
$$
\frac{1}{b_n}\sum_{j=nm+1}^{n(m+1)}
\Lambda\Big( \frac{b_n}{a_n}\sum_{i=1}^{k}\lambda_i
\phi_{l+[nt_{i-1}],[nt_i]-[nt_{i-1}]} \Big)\to
$$
$$
\int\limits_m^{m+1} \zeta\left(
\frac{(1-\alpha)\sum_{i=1}^k\lambda_i\int\limits_{x+t_{i-1}}^{x+t_{i}}p
  y^{-\alpha}dy}{\Bigl|(1-\alpha)\sum_{i=1}^k\lambda_i
\int\limits_{x+t_{i-1}}^{x+t_{i}}p  y^{-\alpha}dy \Bigr|}\right)
\Bigl|(1-\alpha)\sum_{i=1}^k\lambda_i
\int\limits_{x+t_{i-1}}^{x+t_{i}}p  y^{-\alpha}dy \Bigr|^\beta
$$
(if the argument of the function $\zeta$ is $0/0$, then the integrand
is set to be equal to zero), and we treat the other ranges of $m$ in a
manner similar to what has been done in part (ii). This gives us the
stated limit. For $\alpha=1$ we have for any 
$m=1,\ldots, k$, by the regular variation at infinity, 
$$
\frac{1}{b_n} \sum_{j=-[nt_m]+1}^{[nt_{m-1}]}
\Lambda\Big( \frac{b_n}{a_n}\sum_{i=1}^{k}\lambda_i
\phi_{j+[nt_{i-1}],[nt_i]-[nt_{i-1}]}\Big) \to \int\limits_{-t_m}^{-t_{m-1}}
\zeta\left( \frac{ \lambda_m}{|\lambda_m|}\right) |\lambda_m|^\beta\, dx\,,
$$
and so 
$$
\frac{1}{b_n}\sum_{j= -n+1}^{0} \Lambda\Big(
\frac{b_n}{a_n}\sum_{i=1}^{k}\lambda_i
\phi_{j+[nt_{i-1}],[nt_i]-[nt_{i-1}]}\Big) \to \sum_{m=1}^k
(t_m-t_{m-1}) \zeta\left( \frac{ \lambda_m}{|\lambda_m|}\right)
|\lambda_m|^\beta\,, 
$$
while the sum over the rest of the range of $j$ contributes only terms
of a smaller order. Hence the result. 
\end{proof}

\begin{remark}\label{1dimldpmodr}
\emph{ The argument in the proof shows also that the statements of all
  three parts of the lemma remain true if the sums
  $\sum_{l=-\infty}^\infty$ are replaced by sums $\sum_{l=-A_n}^{A_n}$
  with $n/A_n\to 0$ as $n\to\infty$.} 
\end{remark}

\begin{lemma} \label{l:funct.rate.long}
For $1/2<\alpha<1$, let $h_{t_1,\ldots,t_k}$ be defined by \eqref{e:h},
and $\Lambda^{rl}_{t_1,\cdots,t_k}$ defined by
\eqref{e:lambda.rl}. Then for any function $f$ of bounded variation on
$[0,1]$ satisfying $f(0)=0$, 
$$
\sup_{j\in J}  (\Lambda^{rl}_{t_1,\ldots,t_{|j|}})^*\big(f(t_1),f(t_2)-f(t_1),
\ldots,f(t_{|j|})-f(t_{|j|-1})\big)
$$
$$ = \left\{ \begin{array}{lcl}
\Lambda_\alpha^\ast(f^\prime) & if & f\in \mathcal{AC}, \\ \infty &
 & otherwise,\\ \end{array}\right.
$$
where $\Lambda_\alpha^\ast$ is defined by \eqref{e:transf.alpha}. 
\end{lemma}
\begin{proof}
First assume that $f\in \mathcal{AC}$. It is easy to see that the
inequality $\Lambda^*_\alpha(f^\prime) \ge \sup_{j\in
J}(\Lambda^{rl}_{t_1,\cdots,t_{|j|}})^*(f(t_1),f(t_2)-f(t_1),
\ldots,f(t_{|j|}))$ holds 
by considering a function $\psi\in L_\infty[0,1]$, which takes the
value $\lambda_i$ in the interval $(t_{i-1},t_i]$. For the other
inequality, we start by observing that the supremum in the definition
of $\Lambda_\alpha^\ast$ in \eqref{e:transf.alpha} is achieved over
those $\psi\in L_\infty[0,1]$ for which  
the integral 
$$
I_x = \int_0^1 \psi(t)(1-\alpha)|x+t|^{-\alpha}\Bigl[ 
pI_{[x+t\ge 0]}+qI_{[x+t<0]}\Bigr]\, dt \in \mathcal{F}_\Lambda
$$
for almost all real $x$, 
and, hence, also over those $\psi\in L_\infty[0,1]$ for whch $I_x\in
\mathcal{F}_\Lambda^\circ$ for almost every $x$. 

For any $\psi$ as above choose  a sequence of uniformly bounded 
functions $\psi^n$ converging to $\psi$ almost everywhere  on
$[0,1]$, such that for every $n$, $\psi^n$ is of the form
$\sum_i\lambda_i^{n}I_{A_i^{n}},$ where
$A_i^{n}=(t_{i-1}^{n},t_i^{n}]$, for some
$0<t_1^{n}<t_2^{n}<\cdots<t_{k_n}^{n}=1$. Then by the continuity of
$\Lambda$ over $\mathcal{F}_\Lambda^\circ$ and Fatou's lemma, 
\[ \int_0^1\psi(t)f^\prime(t) dt - \int_{-\infty}^\infty 
\Lambda\left( \int_0^1 \psi(t)(1-\alpha)|x+t|^{-\alpha}\Bigl[ 
pI_{[x+t\ge 0]}+qI_{[x+t<0]}\Bigr]\, dt\right)dx\]
\begin{eqnarray*}
&= &\int_0^1\lim_n\psi^n(t)f^\prime(t) dt\\
& & -  \int_{-\infty}^\infty 
\Lambda\left( \int_0^1 \lim_n\psi^n(t)(1-\alpha)|x+t|^{-\alpha}\Bigl[ 
pI_{[x+t\ge 0]}+qI_{[x+t<0]}\Bigr]\, dt\right)dx
\end{eqnarray*}
\begin{eqnarray*}
&= & \lim_n\int_0^1\psi^n(t)f^\prime(t) dt\\
& & - \int_{-\infty}^\infty\lim_n 
\Lambda\left( \int_0^1 \psi^n(t)(1-\alpha)|x+t|^{-\alpha}\Bigl[ 
pI_{[x+t\ge 0]}+qI_{[x+t<0]}\Bigr]\, dt\right)dx
\end{eqnarray*}
\begin{eqnarray*}
&\le  & \lim_n\int_0^1\psi^n(t)f^\prime(t) dt\\
& & -\limsup_n \int_{-\infty}^\infty 
\Lambda\left( \int_0^1 \psi^n(t)(1-\alpha)|x+t|^{-\alpha}\Bigl[ 
pI_{[x+t\ge 0]}+qI_{[x+t<0]}\Bigr]\, dt\right)dx\\
&=& \liminf_n\left\{ \sum_{i=1}^{k_n}\lambda_i^{n}\cdot \big( f(t_i^{n})-f(t_{i-1}^{n})\big)-\Lambda^{rl}_{t_1^{n}, \cdots,t_n^{n}} (\lambda_1^{n},\cdots,\lambda_n^{n})\right\}\\
& \le & \sup_{j\in J}(\Lambda^{rl}_{t_1,\cdots,t_{|j|}})^*\big(f(t_1),f(t_2)-f(t_1), \ldots,f(t_{|j|})-f(t_{|j|-1})\big).
\end{eqnarray*}

Now suppose that $f$ is not absolutely continuous. That is, there
exists $\epsilon>0$  and $0\le r_1^n<s_1^n\le r_2^n<\cdots \le
r_{k_n}^n<s^n_{k_n}\le 1$, such that
$\sum_{i=1}^{k_n}(s_i^n-r_i^n)\rightarrow 0$ but
$\sum_{i=1}^{k_n}|f(s_i^n)-f(r_i^n)|\ge \epsilon$. Let $j^n$ be such
that $t^n_{2p}=s^n_p$ and $t^n_{2p-1}=r^n_p$ (so that
$|j^n|=2k_n$). Now   
\begin{eqnarray*}
&&  \sup_{j\in
J}(\Lambda^{rl}_{t_1,\cdots,t_{|j|}})^*\big(f(t_1),f(t_2)-f(t_1),
\ldots,f(t_{|j|})-f(t_{|j|-1})\big)\\ 
&\ge & \limsup_n \left\{\sup_{\underline{\lambda}^n\in
\mathbb{R}^{2k_n}} \sum_{i=1}^{2k_n}\lambda_i^n\cdot
\big(f(t_i^n)-f(t_{i-1}^n)\big) -
\Lambda^{rl}_{t_1,\cdots,t_{2k_n}}(\underline{\lambda}^n) \right\}\\ 
& \ge & \limsup_n \left\{ A
\sum_{i=1}^{k_n}\big|f(s_i^n)-f(r_i^n)\big|-\Lambda^{rl}_{t_1,
\cdots,t_{2k_n}}(\underline{\lambda}^{n*})\right\}\ge A\epsilon, 
\end{eqnarray*}
where $\lambda^{n*}_{2p-1}=0$ and
$\lambda^{n*}_{2p}=A\big(f(s_i^n)-f(r_i^n)\big)/ |f(s_i^n)-f(r_i^n)|$
($=0$ if $f(s_i^n)-f(r_i^n)=0$) for any $A>0$. The last inequality
follows from an application of dominated convergence theorem,
quadratic behaviour of $\Lambda$ at $0$ and the fact that
$h_{t_1,\cdots,t_{2k_n}}(x;\underline{\lambda}^{n*})\rightarrow 0$ as
$n\rightarrow \infty$ for every $x\in \mathbb{R}$. This completes the
proof since $A$ is arbitrary.  
\end{proof}

\bibliographystyle{elsart-harv}
\bibliography{/Users/souvik/Documents/Miscellanea/Latex/bibfile}

\end{document}